\documentclass[final,onefignum,onetabnum]{siamart220329}

\usepackage[utf8]{inputenc}
\usepackage[algo2e]{algorithm2e}
\usepackage{amssymb}
\usepackage{amsfonts}
\usepackage{amssymb}
\usepackage{caption}
\usepackage{dsfont}
\usepackage{graphicx}
\usepackage{makecell}
\usepackage{optidef}
\usepackage{rotating}
\usepackage{subcaption}
\usepackage{tabularx}
\usepackage{url}
\usepackage{wrapfig}
\usepackage[a4paper,
            bindingoffset=0.2in,
            left=1in,
            right=1in,
            top=1in,
            bottom=1in,
            footskip=.25in]{geometry}
\usepackage[overload]{textcase}

\newcommand{\norm}[1]{\left\lVert#1\right\rVert} 
\newcommand{\quadmat}[3]{{\scriptsize\begin{pmatrix} {#1} & {#2} \\ {#2}^\top & {#3} \end{pmatrix}}}
\newcommand{\quadmatrow}[3]{{\scriptsize\begin{pmatrix} #1 & #2^\top \\ #2 & #3 \end{pmatrix}}}
\newcommand{\Quadmat}[3]{{\scriptsize\begin{pmatrix} #1 & #2 \\ #2 & #3 \end{pmatrix}}}
\newcommand{\innerMat}[2]{\left \langle #1 , \, #2 \right \rangle}
\RestyleAlgo{ruled}

\newcommand{\mbqpNvars}{{n}}
\newcommand{\mbqpNconstrs}{{m}}
\newcommand{\cpDims}{{m}}
\newcommand{\shpVec}{{a}}
\newcommand{\shpConst}{{b}}
\newcommand{\ndisc}{{K}}
\newcommand{\ndiscbin}{{k}}
\newcommand{\Dmat}{{\mathcal{D}}}
\newcommand{\basis}[1]{{\mathds{1}_{\{#1\}}}}

\usepackage{bm,color}
\usepackage{hyperref}
\hypersetup{
    colorlinks,
    linkcolor={red!50!black},
    citecolor={blue!50!black},
    urlcolor={blue!80!black}
}
\usepackage{tikz}
\definecolor{lime}{HTML}{A6CE39}
\DeclareRobustCommand{\orcidicon}{
	\begin{tikzpicture}
	\draw[lime, fill=lime] (0,0) 
	circle [radius=0.16] 
	node[white] {{\fontfamily{qag}\selectfont \tiny ID}};
	\draw[white, fill=white] (-0.0625,0.095) 
	circle [radius=0.007];
	\end{tikzpicture}
	\hspace{-2mm}
}
\foreach \x in {A, ..., Z}{\expandafter\xdef\csname orcid\x\endcsname{\noexpand\href{https://orcid.org/\csname orcidauthor\x\endcsname}
			{\noexpand\orcidicon}}
}

\newcommand\blfootnote[1]{%
  \begingroup
  \renewcommand\thefootnote{}\footnote{#1}%
  \addtocounter{footnote}{-1}%
  \endgroup
}

\newcommand{\SUaffil}{$*$}
\newcommand{\USRAaffil}{$\dagger$}
\newcommand{\QUAILaffil}{$\ddagger$}
\newcommand{\PUaffil}{$\S$}

\makeatletter
\newcommand{\specificthanks}[1]{\@fnsymbol{#1}}%
\makeatother

\title{A copositive framework for \\analysis of hybrid Ising-classical algorithms}
\author{Robin Brown\textsuperscript{\SUaffil \USRAaffil \QUAILaffil}\and David E. Bernal Neira\textsuperscript{\USRAaffil \QUAILaffil \PUaffil}
\and Davide Venturelli\textsuperscript{\USRAaffil \QUAILaffil}
\and Marco Pavone\textsuperscript{\SUaffil}}

\date{\today}

\begin{document}
\maketitle
\blfootnote{\textsuperscript{\SUaffil} Stanford University, Autonomous Systems Laboratory}
\blfootnote{\textsuperscript{\USRAaffil} USRA Research Institute for Advanced Computer Science (RIACS)}
\blfootnote{\textsuperscript{\QUAILaffil} NASA Quantum Artificial Intelligence Laboratory (QuAIL)}
\blfootnote{\textsuperscript{\PUaffil} Purdue University, Davidson School of Chemical Engineering}
\begin{abstract}
    Recent years have seen significant advances in quantum/quantum-inspired technologies capable of approximately searching for the ground state of Ising spin Hamiltonians. 
    The promise of leveraging such technologies to accelerate the solution of difficult optimization problems has spurred an increased interest in exploring methods to integrate Ising problems as part of their solution process, with existing approaches ranging from direct transcription to hybrid quantum-classical approaches rooted in existing optimization algorithms.
    While it is widely acknowledged that quantum computers should augment classical computers, rather than replace them entirely, comparatively little attention has been directed toward deriving analytical characterizations of their interactions.
    In this paper, we present a formal analysis of hybrid algorithms in the context of solving mixed-binary quadratic programs (MBQP) via Ising solvers.
    By leveraging an existing completely positive reformulation of MBQPs, as well as a new strong-duality result, we show the exactness of the dual problem over the cone of copositive matrices, thus allowing the resulting reformulation to inherit the straightforward analysis of convex optimization.
    We propose to solve this reformulation with a hybrid quantum-classical cutting-plane algorithm.
    Using existing complexity results for convex cutting-plane algorithms, we deduce that the classical portion of this hybrid framework is guaranteed to be polynomial time.
    This suggests that when applied to NP-hard problems, the complexity of the solution is shifted onto the subroutine handled by the Ising solver.
\end{abstract}

\section{Introduction}

Recent years have seen significant advances in quantum and quantum-inspired Ising solvers, such as quantum annealers~\cite{johnson2011quantum}, quantum approximate optimization circuits~\cite{farhi2014quantum}, or coherent Ising machines~\cite{honjo2021100}.
These are devices/methodologies designed to heuristically compute solutions of optimization problems of the form: $\min_{z \in \{-1, 1\}^n} \sum_{i, j} J_{i, j} z_i z_j + \sum_{i} h_{i} z_i$, where $J_{i, j}, \, h_i$ are real coefficients and $z_i \in \{-1, 1\}$ are discrete variables to be optimized over.
The promise of leveraging such technologies to speed up the solution of complex optimization problems has spurred many researchers to explore how Ising solvers can be applied to problems in various domains.

A standard approach has emerged where an optimization problem is directly transcribed into an Ising problem, and the returned solution is taken at face value or with minimal post-processing.
While this method works well for problems that organically have an Ising form, sequences of reformulations can result in ill-conditioning of the problem in terms of the number of additional variables, the coupling strengths, and the optimization landscape. 
Most unnaturally, however, relying solely on the Ising solver means forgoing the advantages of already powerful classical computers.

In an effort to introduce meaningful interplay between classical and quantum computers, a few authors have proposed decomposition methods based on the Alternating Direction Method of Multipliers (ADMM), \cite{gambella2020multiblock}, or Benders Decomposition (BD), \cite{chang2020hybrid, zhao2021hybrid,  paterakis2021hybrid}.
Critically, when ADMM is applied to non-convex problems, it is not guaranteed to converge.
When it does, it is often to a local optimum without convergence guarantees to global optimality.
On the other hand, while it may be possible to derive optimality guarantees using BD, proving convergence typically relies on an exhaustive search through the ``complicating variables".
This makes it unclear whether such an algorithmic scaffold is primed to take advantage of speed-ups that the Ising solver may offer.
\paragraph{Contributions}
Our work is motivated by the desire to rigorously analyze the interplay between classical and quantum machines in hybrid algorithms.
Such analysis is a cornerstone for articulating and setting standards for hybrid quantum-classical optimization algorithms.
Specifically, we espouse convergent hybrid quantum-classical algorithms that (1) use Ising solvers as a primitive while offering some resilience to their heuristic nature and (2) have polynomial complexity in the classical portions of the algorithm. 
To this end, the contribution of this paper is an algorithmic framework that satisfies these key desiderata.
Concretely, 
\begin{enumerate}
    \item We revisit and prove strong duality of a result in \cite{Burer2009} to show that the convex copositive formulation of many mixed-binary quadratic optimization problems is exact.
    Neglecting the challenges of working with copositive matrices, convex programs are a well-understood class of optimization problems with a wide variety of efficient solution algorithms.
    By reformulating mixed-binary quadratic programs as copositive programs, we open the door for hybrid-quantum classical algorithms that are based on existing convex optimization algorithms.
    \item To solve the copositive programs, we propose applying a standard cutting-plane algorithm, which we modify in a novel way using a hybrid quantum-classical approach.
    In particular, the cutting-plane algorithm serves as a template for a solution algorithm that alternates between checking copositivity and other operations.
    We hybridize the algorithm by approximating the copositivity checks via discretization and solving them with an Ising solver.
    We show that the complexity of the portion of the algorithm handled by the classical computer has polynomial scaling.
    This analysis suggests that when applied to NP-hard problems, the complexity of the solution is shifted onto the subroutine handled by the Ising solver.
    \item We conducted benchmarking based on the maximum clique problem to validate our theoretical claims and evaluate potential speed-ups from using a stochastic Ising solver in lieu of a state-of-the-art deterministic solver or an Ising heuristic. Results indicate that the Ising formulation of the subproblems of the hybrid algorithm is efficient versus a MIP formulation in \texttt{Gurobi}, and the hybrid algorithm potentially is competitive even against a non-hybridized Ising formulation of the full problem solved by simulated annealing.
\end{enumerate}
We emphasize that the contribution of this work is not the copositive reformulation or a novel cutting-plane algorithm but rather the insight that these ideas are synergistic with recent advances in quantum(-inspired) computing.
In particular, copositive optimization is useful for deriving and analyzing new hybrid algorithms rooted in existing convex optimization algorithms, thus filling a gap in the hybrid algorithms literature.

While preparing this manuscript, a hybrid quantum-classical method relying upon a Frank-Wolfe method was published \cite{yurtsever2022q}. This work also leverages a similar copositive reformulation of quadratic binary optimization problems.
We highlight the differences between that work and the one in our manuscript below.
This manuscript considers the optimization problem class of mixed-binary quadratic programs, while in \cite{yurtsever2022q}, the authors propose their method for quadratic binary optimization problems, a subcase of the problems considered herein.
Moreover, we provide proof of the exactness and strong duality for copositive/completely positive optimization stemming from the mixed-integer quadratic reformulation, addressing an open question in the field.
In their manuscript, \cite{yurtsever2022q} conjectures the results proved in this manuscript to be true.
Finally, our solution method, which is based on cutting-plane algorithms, has a potential exponential speed-up in runtime compared to Frank-Wolfe algorithms.

\subsection{Related Work}
One dominant method for mapping optimization problems into Ising problems is through direct transcription.
This process typically involves discretizing continuous variables and passing constraints into the objective through a penalty function; the returned solution is often taken at face value or with minimal post-processing to enforce feasibility.
Owing to its simplicity, this process has found applications in a variety of problems, including jobshop scheduling \cite{venturelli2016job}, routing problems \cite{harwood2021formulating}, community detection \cite{negre2020detecting}, and all of Karp's 21 NP-complete models \cite{lucas2014ising}, among others.
Critically, unless a problem organically takes an Ising form, this approach often requires many auxiliary variables (spins), introduces large skews in the coupling coefficients, and can result in poor conditioning of the optimization landscape, thus limiting the problems that can be solved on near-term devices.
Similarly, extending the class of applicable problem instances requires deriving increasingly complex sequences of reformulations, each of which reduces the solubility of the final reformulation. 
More importantly, an algorithm with minimal interactions between classical and quantum computers disregards the bountiful successes of classical computers in the past decades.
This has inspired some researchers to examine how quantum computers can be used to augment classical computers instead of replacing them entirely \cite{callison2022hybrid}. 

As an alternative to direct transcription, there is a burgeoning body of literature exploring the potential of decomposition methods for designing hybrid quantum-classical algorithms.
These generally refer to algorithms that divide effort between a classical and quantum computer, with each computer informing the computation carried out by the other.
Among these, algorithms based on the Benders Decomposition (BD) are gaining traction.
BD is particularly effective for problems characterized by ``complicating variables", for which the problem becomes easy once these variables are fixed. 
For example, a mixed-integer linear program (MILP) becomes a linear program (LP) once the integer variables are fixed--the integers are the complicating variables. 
BD iterates between solving a master problem over the complicating variables and sub-problems where the complicating variables are fixed, whose solution is used to generate cuts for the master problem.
Both \cite{chang2020hybrid} and \cite{zhao2021hybrid} consider mixed-integer programming (MIP) problems where the integer variables are linked to the continuous variables through a polyhedral constraint and leverage a reformulation where dependence on the continuous variables is expressed as constraints over the extreme rays and points of the original feasible region.
Because the number of extreme rays and points may be exponentially large, the constraints are not written down in full but are iteratively generated from the solutions of the sub-problems.
The master problem is an integer program consisting of these constraints and is solved using the quantum computer.
Notably, the generated constraint set may be large, with the worst case being the generation of the entire constraint set, resulting in a large number of iterations.
The approach in \cite{paterakis2021hybrid} attempts to mitigate this by generating multiple cuts per iteration and selecting the most informative subset of these cuts. 
Instead of using the quantum computer to solve the master problem, the quantum computer is used to heuristically select cuts based on a minimum set cover or maximum coverage metric.
While this may effectively reduce the number of iterations and size of the constraint set, the master problem is often an integer program that may be computationally intractable.
For each of the proposed approaches, it is unclear how the complexity of the problem is distributed through the solution process
For example, for \cite{chang2020hybrid, zhao2021hybrid} the complexity might show up in the number of iterations, and for \cite{paterakis2021hybrid} it might show up when solving the master problem.
Consequently, it is ambiguous whether BD-based approaches can take advantage of a speed-up in the Ising solver, even if one were to exist.

Another decomposition that has been explored is based on the Alternating Direction Method of Multipliers (ADMM)\cite{gambella2020multiblock}.
This is an algorithm to decompose large-scale optimization problems into smaller, more manageable sub-problems \cite{boyd2011distributed}.
While originally designed for convex optimization, ADMM has shown great success as a heuristic for non-convex optimization as well, \cite{diamond2018general}, and significant progress has been made towards explaining its success in such settings \cite{wang2019global}.
In \cite{gambella2020multiblock}, the authors propose an ADMM-based decomposition with three sub-problems: the first being over just the binary variables, the second being the full problem with a relaxed copy of the binary variables, and the third being a term that ties the binary variables and their relaxed copies together.
For quadratic pure-binary problems, the authors show that the algorithm converges to a stationary point of the augmented Lagrangian, which may not be a global optimizer
Convergence to a global optimum is only guaranteed under the more stringent Kurdyka-\L ojasiewicz conditions on the objective function \cite{attouch2013convergence}.
Unfortunately, the assumptions guaranteeing convergence to a stationary point fail in the presence of continuous variables. 

A third class of decomposition proposed and implemented in the \texttt{qbsolv} solver is based on tabu search \cite{BoothReinhardtEtAl2017}.
\texttt{qbsolv} can be seen as iterating between a large-neighborhood local search (using an Ising solver) and tabu improvements to locally refine the solution (using a classical computer), where previously found solutions are removed from the search space in each iteration.
During the local search phase, subsets of the variables are jointly optimized while the remaining variables are fixed to their current values.
The solution found in this phase is then used to initialize the tabu search algorithm, and the process is repeated for a fixed number of iterations.
Critically, it is unclear whether the algorithm is guaranteed to converge and, if so, what its optimality guarantees are.
While finite convergence of tabu search is investigated in \cite{glover2002tabu}, it relies on either recency or frequency memory that ensures an exhaustive search of all potential solutions.

Another approach for purely integer programming problems is based on the computation of a Graver basis through the computation of the integer null-space of the constraint set as proposed in \cite{alghassi2019graver}.
This null-space computation is posed as a quadratic unconstrained binary optimization (QUBO) and then post-processed to obtain the Graver basis of the constraint set, a test-set of the problem.
The test-set provides search directions for an augmentation-based algorithm. 
For a convex objective, it provides a polynomial oracle complexity in converging to the optimal solution.
The authors initialize the problem by solving a feasibility-based QUBO and extend this method to non-convex objectives by allowing multiple starting points for the augmentation.
The multistart procedure also alleviates the requirement for computing the complete Graver basis of the problem, which grows exponentially with the problem's size.
Considering an incomplete basis or non-convex objectives makes the Graver Augmentation Multistart Algorithm (GAMA) a heuristic for general integer programming problems, and it cannot address problems with continuous variables.

In this paper, we seek to address a gap in the literature on a rigorous theory of hybrid quantum-classical optimization. 
By revisiting the hidden convex structure of non-convex problems, we pave the way for hybrid algorithms based on efficient convex optimization.
We show that algorithms derived through this approach inherit the straightforward analysis of convex optimization without sacrificing the potential benefits of quantum computing for non-convex problems.

While there is optimism regarding improvements to and our understanding of quantum technology in the coming decades, few expect that they will replace classical computers entirely.
We believe that the method presented in this paper is an approach to algorithm design that anticipates a future where quantum computers and classical computers work in tandem.
In particular, we envision a mature theory of hybrid algorithms that clearly delineates how quantum and classical computers should complement each other.

\subsection{Quantum/Quantum-inspired Ising Solvers}
Adiabatic quantum computing (AQC) is a quantum computation paradigm that operates by initializing a system in the ground state of an initial Hamiltonian (i.e., the optimal solution of the corresponding objective function) and slowly sweeping the system to an objective Hamiltonian.
This Hamiltonian, referred to as the cost Hamiltonian, maps the objective function of the classical Ising model onto a system with as many quantum bits, or qubits, as original variables in the Ising model.
The adiabatic theorem of quantum mechanics states that if the system evolution is ``sufficiently slow", the system ends up in the ground state of the desired Hamiltonian.
Here, ``sufficiently slow" depends on the minimum energy gap between the ground and the first excited state throughout the system evolution \cite{albash2018adiabatic}.
Since the evaluation of the minimal gap is mostly intractable, one is forced to phenomenologically ``guess'' the evolution’s speed, and if it is too fast, the undesired non-adiabatic transitions can occur.
Additionally, real devices are plagued with various incarnations of physical noise, such as thermal fluctuations or decoherence effects, that can hamper computation.
The situation is further exacerbated by the challenge of achieving dense connectivity between qubits.
Densely connected problems are embedded in devices by chaining together multiple physical qubits to represent one logical qubit. 
The heuristic computational paradigm that encompasses the additional noise and non-quantum effects is known as Quantum Annealing (QA).
\cite{hauke2020perspectives} provides a review on QA with a focus on possible routes towards solving the open questions in the field. 

An alternative paradigm to AQC is the gate-based model of quantum computing.
Within the gate-based model, Variational Quantum Algorithms (VQAs) is a class of hybrid quantum-classical algorithms that can be applied to optimization \cite{cerezo2021variational}.
VQAs share a common operational principle where the ``loss function" of a parameterized quantum circuit is measured on a quantum device and evaluated on a classical processor, and a classical optimizer is used to update (or ``train") the circuit's parameters to minimize the loss.
VQAs are often interpreted as a quantum analog to machine learning, leaving many similar questions open regarding their trainability, accuracy, and efficiency.
Most similar in spirit to this work is a theory of variational hybrid quantum-classical algorithms proposed in \cite{mcclean2016theory}.
However, they primarily focus on algorithmic improvements to the quantum portion, with discussion of the classical optimization being limited to empirical evaluations of existing derivative-free optimization algorithms.
More recently, \cite{bittel2021training} analyzed the complexity of training VQAs, and through reductions from the maximum cut problem, showed that it is NP-hard.
The analysis presented in this paper complements these prior works in developing a more complete picture of the interplay between quantum and classical computers.

The quantum approximate optimization algorithm (QAOA) is a specific instance of a VQA where the structure of the quantum circuit is the digital analog of adiabatic quantum computing \cite{farhi2014quantum}. 
QAOA operates by alternating the application of the cost Hamiltonian and a mixing Hamiltonian; the number of alternating blocks is referred to as the circuit depth.
For each one of the alternating steps, either mixing or cost application, a classical optimizer needs to determine how long each step should be performed, encoded as rotation angles.
Optimizing the expected cost function with respect to the rotation angles is a continuous low-dimensional non-convex problem.
QAOA is designed to optimize cost Hamiltonians, such as the ones derived from classical Ising problems.
Performance guarantees can be derived for QAOA with well-structured problems, given that the optimal angles are found in the classical optimization step.
Although approximation guarantees have not been derived for arbitrary cost Hamiltonians, even depth-one QAOA circuits have non-trivial performance guarantees for specific problems and cannot be efficiently simulated on classical computers  \cite{farhi2016quantum}, thus bolstering the hope for a speed-up in near-term quantum machines.
Moreover, the algorithm's characteristics, such as relatively shallow circuits, make it amenable to be implemented in currently available noisy intermediate-scale quantum (NISQ) computers compared to other algorithms requiring fault-tolerant quantum devices~\cite{preskill2018quantum}.
While QAOA's convergence to optimal solutions is known to improve with increased circuit depth and to succeed in the infinite depth limit following its equivalence to AQC, its finite depth behavior has remained elusive due to the challenges in analyzing quantum many-body dynamics and other practical complications such as decoherence when implementing long quantum circuits, compilation issues, and hardness of the optimal angle classical problem~\cite{uvarov2021barren}.
Even considering these complications, QAOA has been extensively studied and implemented in current devices~\cite{willsch2020benchmarking,harrigan2021quantum}, becoming one of the most popular alternatives to address combinatorial optimization problems modeled as Ising problems using gate-based quantum computers.
Several other quantum heuristics for Ising problems have been proposed, usually requiring fault-tolerant quantum computers.
We direct the interested reader to a recent review on the topic~\cite{sanders2020compilation}.

An alternative physical system for solving Ising problems that has emerged is coherent Ising machines (CIMs), which are optically pumped networks of coupled degenerate optical parametric oscillators.
As the pump strength increases, the equilibrium states of an ideal CIM correspond to the Ising Hamiltonian's ground states encoded by the coupling coefficients. 
Large-scale prototypes of CIMs have achieved impressive performance in the lab, thus driving the theoretical study of their fundamental operating principles. 
While significant advances have been made on this front, we still lack a clear theoretical understanding of the CIMs' computational performance. Since a thorough understanding of the CIM is limited by our capacity to prove theorems about complex dynamic systems, near-term usage of CIMs must treat them as a heuristic rather than a device with performance guarantees \cite{yamamoto2017coherent}.
Even so, there are empirical observations that in many cases, the median complexity of solving Ising problems using CIM scales as $\exp{\sqrt{N}}$ where $N$ is the size of the problem~\cite{mohseni2022ising}, making it a potential approach to solve these problems efficiently in practice.
We note that there are other types of Ising machines, including classical thermal annealers (based on magnetic devices \cite{shim2017ising}, optics \cite{pierangeli2019large}, memristors \cite{bojnordi2016memristive}, and digital hardware accelerators \cite{matsubara2020digital}), dynamical-systems solvers (based on optics \cite{mcmahon2016fully} and electronics \cite{chou2019analog}), superconducting-circuit quantum annealers \cite{albash2018demonstration}, and neutral atoms arrays \cite{henriet2020quantum}.
We direct the interested reader to \cite{mohseni2022ising}, which provides a recent review and comparison of various methods for constructing Ising machines and their operating principles.

\paragraph{Organization}
In Section \ref{sec:prelim}, we present notation, terminology, and the problem setting covered by our approach.
In Section \ref{sec:approach}, we introduce the proposed framework, including convex reformulation via copositive programming, a high-level overview of cutting-plane algorithms, and a specific discussion of their application to copositive programming.
Section \ref{sec:experiments} provides numerical experiments supporting our assertions about the proposed approach.
Finally, we conclude and highlight future directions in Section \ref{sec:conclusion}.

\section{Preliminaries}\label{sec:prelim}
\subsection{Notation and Terminology}
In this paper, we solely work with vectors and matrices defined over the real numbers and reserve lowercase letters for vectors and uppercase letters for matrices. 
We will also follow the convention that a vector $x \in \mathbb{R}^n$ is to be treated as a column vector, i.e., equivalent to a matrix of dimension $n \times 1$.
For a matrix $M$, we use $M_{i, j}$ to denote the entry in the $i$th row and $j$th column, $M_{i, *}$ denotes the entire $i$th row, and $M_{*, j}$ denotes the entire $j$th column.
In the text, we frequently use block matrices with structured zero entries; we use $\cdot$ as shorthand for zero entries.
We use $\mathds{1}$ to denote the all-ones vectors and $\basis{j}$ to denote the $j$th standard basis vector (i.e., a vector where all entries are zero except for a 1 for the jth entry).
The $p$-norm of a vector $v \in \mathbb{R}^n$ is defined as $\norm{v}_p := \left(\sum_{i= 1}^n v_i^p \right)^{1/p}$.
We reserve the letter $I$ to denote the identity matrix.
For two matrices, $M$ and $N$, we use $\langle M, \, N \rangle = \text{Tr}(M^\top N)$ to denote the matrix inner product.
Note that for two vectors, $\text{Tr}(x^\top y) = x^\top y$ because $x^\top y$ is a scalar, so the matrix inner product is consistent with the standard inner product on vectors.
For sets, $S_M + S_N : = \{ M + N \mid M \in S_M, N \in S_N\}$ is their Minkowski sum,  $S_M \cup S_N$ their union, and $S_M \cap S_N$ their intersection.
For a cone, $\mathcal{K}$, its dual cone is defined as $\mathcal{K}^* = \{X \mid \langle X, \, K \rangle \geq 0, \forall K \in \mathcal{K}\}$.
While we work with matrix cones in this paper, this definition of dual cones is consistent with vector cones as well.
In this paper, the two cones we will work with are the cone of completely positive matrices and the cone of copositive matrices.
The cone of completely positive (CP) matrices, $C^*$, is the set of matrices that have a factorization with entry-wise non-negative entries:
\begin{equation}
    \mathcal{C}^*_{n} := \{X \in \mathbb{R}^{n \times n} \mid X = \sum_k x^{(k)} (x^{(k)})^\top, \quad x^{(k)} \in  \mathbb{R}^n_{\scriptscriptstyle \geq 0} \}
\end{equation}
The cone of copositive matrices, $\mathcal{C}$, is the set of matrices defined by:
\begin{equation}
    \mathcal{C}_{n} := \{X \in \mathbb{R}^{n \times n} \mid v^\top X v \geq 0, \quad \forall v \in  \mathbb{R}^n_{\scriptscriptstyle \geq 0} \}
\end{equation}
As suggested by the notation, the cones of completely positive and copositive matrices are duals of each other.
We use $S^n_{++}$ to denote the cone of positive definite matrices.

In this paper, we will use the terms Ising problem and quadratic unconstrained binary optimization (QUBO) interchangeably.
An Ising problem is an optimization problem of the form: $\min_{z \in \{-1, 1\}^n} \sum_{i, j} J_{i, j} z_i z_j + \sum_{i} h_{i} z_i$, where $J_{i, j}, \, h_i$ are real coefficients and $z_i \in \{-1, 1\}$ are discrete variables to be optimized over.
A QUBO, which is an optimization problem of the form $ \min_{x \in \{0, \, 1\}^n} \sum_{i, j} Q_{i, j} x_i x_j$ can be reformulated as an Ising problem using the change of variable $z = 2x - \mathds{1}$.
This translates to coefficients in the Ising problem $J_{i, j} = \frac{1}{4} Q_{i, j}$, $h_i = \frac{1}{2} \sum_{j} Q_{i, j}$, and a constant offset of $\frac{1}{4}\sum_{i, j} Q_{i, j}$.
\subsection{Problem Setting}
In this paper, we consider mixed-binary quadratic programs (MBQP) of the form:
\begin{mini}
{x \in \mathbb{R}^\mbqpNvars}{x^\top Q x + 2 c^\top x}
{\label{eq:MBQP}\tag{MBQP}}{}
\addConstraint{Ax = b, \quad A \in \mathbb{R}^{\mbqpNconstrs \times \mbqpNvars}, \, b \in \mathbb{R}^{\mbqpNconstrs}}
\addConstraint{x \geq 0}
\addConstraint{x_j \in \{0, 1\}, \quad  j \in B}
\end{mini}
where the set $B \subseteq \{1, \ldots, \mbqpNvars\}$ indexes which of the $n$ variables are binary.
This is a general class of problems that encompasses problems including QUBOs, standard quadratic programming, the maximum stable set problem, and the quadratic assignment problem.
Because mapping to an Ising problem can also be equivalently expressed as a QUBO, many problems tackled with Ising solvers thus far pass through a formulation similar to the form of Problem \eqref{eq:MBQP}.
Using the result in \cite[Sec. 3.2]{Burer2009}, the formulation considered in this paper can be extended to include constraints of the form $x_ix_j = 0$ that force at least one of $x_i$ or $x_j$ to be zero, i.e., complementarity constraints.
For ease of notation, this extension is left out of the present discussion.

\section{Proposed Methodology}\label{sec:approach}
In this section, we will discuss our proposed methodology for solving Problem \eqref{eq:MBQP} given access to Ising solvers.
Our result relies on a convex reformulation of Problem \eqref{eq:MBQP} as a copositive program.
Leveraging convexity, we propose to solve the problem using cutting-plane algorithms.
These belong to a broad class of convex optimization algorithms whose standard components give rise to a natural separation between the role of the Ising solver versus a classical computer.

We first state Burer's exact reformulation of Problem \eqref{eq:MBQP} as a completely positive program and its dual copositive program.
We then show that under mild conditions (i.e., feasibility and boundedness) of the original MBQP, the copositive and completely positive programs exhibit strong duality.
We will then introduce the class of cutting-plane algorithms and summarize the complexity guarantees of several well-known variants.
Finally, we explicitly show how cutting-plane algorithms can be used to solve copositive optimization problems given a copositivity oracle and discuss how to implement a copositivity oracle using an Ising solver.

\subsection{Convex formulation as a copositive program}
In his seminal work, Burer showed that MBQPs can be represented exactly as completely positive programs of the form:
\begin{mini}
    {X \in \mathbb{R}^{\mbqpNvars \times \mbqpNvars}, \, x \in \mathbb{R}^\mbqpNvars}{\innerMat{\quadmat{Q}{c}{\cdot}}{\quadmat{X}{x}{1}}}
    {\label{eq:CPP}\tag{CPP}}{}
    \addConstraint{\innerMat{\quadmatrow{\cdot}{\frac{1}{2} A_{i, *}}{\cdot}}{\quadmat{X}{x}{1}} = b_i, \, i = 1, \ldots, \mbqpNconstrs}
    \addConstraint{\innerMat{\Quadmat{A^\top_{i, *} A_{i, *}}{\cdot}{\cdot}}{\quadmat{X}{x}{1}} = b^2_i, \, i = 1, \ldots, \mbqpNconstrs}
    \addConstraint{\innerMat{\quadmat{-\basis{j} \basis{j} ^\top}{\frac{1}{2} \basis{j} }{\cdot}}{\quadmat{X}{x}{1}} = 0, \,  j \in B}
    \addConstraint{\quadmat{X}{x}{1} \in \mathcal{C}^*_{\mbqpNvars + 1},}
\end{mini}
where exactness means that Problems \eqref{eq:MBQP} and \eqref{eq:CPP} have the same optimal objective and for an optimal solution, $(x^*, X^*)$, of \eqref{eq:CPP}, $x^*$ lies within the convex hull of optimal solutions for \eqref{eq:MBQP} \cite[Theorem 2.6]{Burer2009}.
Similar to semi-definite programming (SDP) relaxations, the completely positive formulation involves lifting the variables in \eqref{eq:MBQP} to a matrix variable representing their first and second-degree monomials, making the objective function and constraints linear.
Unlike SDP relaxations, however, the complete positivity constraint is sufficient for ensuring that the feasible region of \eqref{eq:CPP} is exactly the convex hull of the feasible region of \eqref{eq:MBQP}.
This distinction is what ensures that the optimal value of \eqref{eq:CPP} is exactly that of \eqref{eq:MBQP}, whereas for an SDP relaxation, the optimal solution may lie outside of the convex hull of \eqref{eq:MBQP}, resulting in a lower objective value (i.e., a \emph{relaxation gap}).

Taking the dual of \eqref{eq:CPP} yields a copositive optimization problem of the form \cite[Section 5.9]{boyd2004convex}:
\begin{maxi}
{\mu, \lambda, \gamma}{\gamma + \sum_{i = 1}^{\mbqpNconstrs} \mu^{\text{(lin)}}_i b_i + \mu^{\text{(quad)}}_{i} b_i^2}
{\label{eq:COP}\tag{COP}}{}
\addConstraint{M(\mu, \lambda, \gamma) \in \mathcal{C}_{\mbqpNvars + 1},}
\end{maxi}
where
\begin{equation}\label{def:M}
    \begin{aligned}
        M(\mu, \lambda, \gamma) := &\quadmat{Q}{c}{\cdot} - \sum_{i = 1}^{\mbqpNconstrs} \mu^{\text{(lin)}}_i \quadmatrow{\cdot}{\frac{1}{2} A_{i, *}}{\cdot} - \sum_{i = 1}^{\mbqpNconstrs} \mu^{\text{(quad)}}_i \Quadmat{A^\top_{i, *} A_{i, *}}{\cdot}{\cdot}\\
        & - \sum_{j \in B} \lambda_j \quadmat{-\basis{j}  \basis{j} ^\top}{\frac{1}{2} \basis{j} }{\cdot} - \gamma\Quadmat{\cdot}{\cdot}{1}
    \end{aligned}
\end{equation}
is a parametrized linear combination of the constraint matrices.
The dual copositive program has a linear objective and a single copositivity constraint--this is a convex optimization problem.
While weak duality always holds between an optimization problem and its dual, strong duality is not generally guaranteed.
Showing that strong duality holds is critical for ensuring convergence of specific optimization algorithms and exactness when solving the dual problem as an alternative to solving the primal.
\begin{theorem}[Strong Duality]\label{thm:strong_duality}
If Problem \eqref{eq:MBQP} is feasible with bounded feasible region, then strong duality holds between Problems \eqref{eq:CPP} and \eqref{eq:COP} (i.e., $ \min \eqref{eq:CPP} = \max \eqref{eq:COP}$).
\begin{proof}[Proof Sketch]
Our proof proceeds by first showing strong duality between the alternative representation of \eqref{eq:CPP} (using a homogenized formulation of the equality constraints) and its dual.
By showing that the optimal value of \eqref{eq:COP} is lower-bounded by the optimal value of this homogenized dual problem, we can sandwich the optimal values of Problems \eqref{eq:CPP} and \eqref{eq:COP} by those of a primal-dual pair that has been shown to exhibit strong duality.
The complete proof of this result is provided in Appendix \ref{subsec:strong_dual}.
\end{proof}
\end{theorem}
In prior work, characterization of the duality gap between Problems \eqref{eq:CPP} and \eqref{eq:COP} has remained elusive because the feasible region of Problem \eqref{eq:CPP} never has an interior, thus prohibiting straightforward application of Slater's constraint qualification.
This result is significant because it shows that under mild conditions, the copositive formulation is exact.
This means that the optimal values of Problems \eqref{eq:MBQP} and \eqref{eq:COP} are equivalent, so solving Problem \eqref{eq:COP} is a valid alternative to solving Problem \eqref{eq:MBQP}.
The framework developed in this paper will ultimately produce approximate solutions for Problems \eqref{eq:CPP} and \eqref{eq:COP}, which we anticipate can be used to speed up the solution process of a purely classical solver for \eqref{eq:MBQP}.
For example, the heuristic solutions or cuts used to generate them might be used to warm-start or initialize a purely classical solver.

While Problems \eqref{eq:CPP} and \eqref{eq:COP} are both convex, neither resolve the difficulty of Problem \eqref{eq:MBQP} as even checking complete positivity (resp. copositivity) of a matrix is NP-hard (resp. co-NP-complete) \cite{murty1987some}.
Instead, they should be viewed as ``packaging" the complexity of the problem entirely in the copositivity/complete positivity constraint.
There are a number of classical approaches for (approximately) solving copositive/completely positive programs directly, such as the sum of squares hierarchy \cite{Parrilo2000, lasserre2001global}, feasible descent method in the completely positive
cone, approximations of the copositive cone by a sequence of
polyhedral inner and outer approximations, among others \cite{dur2021conic,burer2012copositive,dur2010copositive}.
In this paper, we will exploit the innate synergy between checking copositivity, which is most naturally posed as a quadratic minimization problem, and solving Ising problems.
This perspective is suggestive of a hybrid quantum-classical approach where the quantum computer is responsible for checking feasibility (i.e., the ``hard part") of the copositive program while the classical computer directs the search towards efficiently reducing the search space.

\subsection{Cutting-Plane/Localization Algorithms}\label{subsec:cp_alg}
Cutting-plane/localization algorithms are convex optimization algorithms that divide labor between checking feasibility, abstracted as a \emph{separation oracle}, and optimization of the objective
\footnote{The term ``cutting-plane algorithm" overloaded in the literature, with one class referring very explicitly to those designed for convex/quasi-convex optimization problems (for a pedagogical reference, we refer the interested reader to \cite{boyd2007localization}) and the second referring more broadly to algorithms that iteratively generate cuts (including algorithms for integer programming and non-convex optimization). In this work, we refer specifically to those designed for convex/quasi-convex optimization.}.
In this section, we provide a high-level overview of each algorithmic step and summarize both the runtime and oracle complexities of several well-known variants; these complexity measures will ultimately correspond to the complexity of the sub-routine handled by the classical computer and the number of calls to the Ising solver, respectively.

While cutting-plane algorithms are often used to solve both constrained and unconstrained optimization problems, they are generally evaluated in terms of their complexity when solving the \emph{feasibility problem}.
\begin{definition}[Feasibility Problem]
For a set of interest $S \subset \mathbb{R}^{\cpDims} $, which can only be accessed through a \underline{separation oracle}, the feasibility problem is concerned with either finding a point in the set $x \in S$ or proving that $S$ does not contain a ball of radius $r$.
\end{definition}

\begin{definition}[Separation Oracle]
A separation oracle for a set $S$, $\texttt{Oracle}_S(\cdot)$ takes as input a point $x \in \mathbb{R}^{\cpDims}$ and either returns \texttt{True} if $x \in S$ or a separating hyperplane if $x \not \in S$.
A separating hyperplane is defined by a vector, $\shpVec \in \mathbb{R}^{\cpDims}$ and scalar $\shpConst \in \mathbb{R}$ such that $\shpVec^\top s \leq \shpConst$ for all $s \in S$ but $\shpVec^\top x \geq \shpConst$.
\end{definition}
The feasibility problem formulation is non-restrictive because these methods can be readily adapted to solving quasi-convex optimization problems with only a simple modification to the separation oracle.
In particular, if the separation oracle indicates feasibility and returns a vector $g\in \mathbb{R}^m$ where any vectors $x, y \in \mathbb{R}^m$ with $f(y) < f(x)$ implies that $g^\top y \geq g^\top x$, this serves as a separating hyperplane for the subset of the feasible region that has a better objective than the test point.
If $f$ is subdifferentiable, any subgradient $g \in \partial f(x)$ satisfies this condition, and for Problem \eqref{eq:COP}, choosing $g$ as the objective's coefficient vector is sufficient.

Although there are many variations of cutting-plane algorithms, at a high level, they follow a standard template that consists of alternating between checking feasibility of a test point, updating an outer approximation of the feasible region, and judiciously selecting the next test point.
This standard template is summarized in Algorithm \ref{alg:cp_meta}.
An overview of the Ellipsoid algorithm is included in the Appendix as a representative example of cutting-plane algorithms, and we direct the interested reader to the references listed in Table \ref{tab:cp} for specific implementation details.
By choosing subsequent test points to be the center of the outer approximation, the algorithm is guaranteed to make consistent progress in reducing the search space (where the metric of progress may also vary across cutting-plane algorithms).
Intuitively, cutting plane algorithms can be considered a high-dimensional analog of binary search.
We note that the requirement $\text{Vol}(S_0) \leq R$ means that the initial set must be bounded.
While this is a standard assumption in the cutting-plane literature, finding such an $S_0$ may be non-trivial.
Procedurally, one may construct a bounded outer approximation using a linear program as in Step 1, Algorithm 1 of \cite{sivaramakrishnan2007properties}.

\begin{algorithm2e}
\caption{Cutting-plane meta-algorithm (feasibility problem)}\label{alg:cp_meta}
\KwIn{$S_0 \subseteq \mathbb{R}^{\cpDims}$ (Initial Set) with $\texttt{Vol}(S_0) \leq R$}
\KwOut{$x \in S$ or \texttt{False} if S does not contain a ball of volume $r$}
$x \gets \texttt{Center}(S_0)$\;
$k \gets 0$\;
\While{\texttt{Oracle}(x) is not \texttt{True} and $\texttt{Vol}(S_k) \geq r$}{
    $S_{k + 1} \gets \texttt{Add\_Cut}(S_k, \texttt{Oracle}(x))$\;
    $x \gets \texttt{Center}(S_{k + 1})$\;
    $k \gets k + 1$\;
}
\uIf{\texttt{Oracle}(x) is \texttt{True}}{
\Return{$x$}\;
}
\Else{
\Return{\texttt{False}}\;
}
\end{algorithm2e}

A number of well-known variants of cutting-plane algorithms are summarized in Table \ref{tab:cp}.
Differences across instantiations of cutting-plane algorithms vary in how subsequent test points are chosen, how the outer approximation is updated, and how progress in decreasing the outer approximation's size is measured.
Each of the surveyed variants strikes a different balance between the computational effort needed to compute a good center versus the resolution used to represent the outer approximation.
Critically, except for the Center of Gravity method, all cutting-plane algorithms summarized in Table \ref{tab:cp} have a polynomial complexity in the dimension of the optimization variables in terms of both oracle queries and total runtime excluding the oracle calls (i.e., the total complexity of adding the cuts and generating test points).
This suggests that if a cutting-plane algorithm were applied to Problem \eqref{eq:COP}, the complexity of the problem is offloaded onto the separation oracle; this is the subroutine we propose to handle using an Ising solver.
\begin{table}[h!]
    {
    \centering
    \begin{tabular}{|c|c|c|c|}
    \hline
    Name & Oracle Queries & \makecell{Total runtime\\ (excluding oracle queries)} & References\\
    \hline
    Center of Gravity & $\mathcal{O}({\cpDims} \log(\frac{R}{r}))$ & \#P-hard \cite{rademacher2007approximating} & \cite{levin1965algorithm}\\
    Ellipsoid & $\mathcal{O}({\cpDims}^2 \log({\cpDims}\frac{R}{r}))$ & $\mathcal{O}({\cpDims}^4 \log({\cpDims}\frac{R}{r}))$ & \cite{shor1977cut, yudin1976evaluation,  khachiyan1980polynomial}\\
    Inscribed Ellipsoid & $\mathcal{O}({\cpDims} \log ({\cpDims} \frac{R}{r}))$ & $\mathcal{O}(({\cpDims} \log ({\cpDims} \frac{R}{r}))^{4.5})$ & \cite{khachiyan1988method, nesterov1989self} \\ 
    Volumetric Center & $\mathcal{O}({\cpDims} \log ({\cpDims} \frac{R}{r}))$ & $\mathcal{O}({\cpDims}^{1 + \omega} \log ({\cpDims} \frac{R}{r}))$ & \cite{vaidya1989new} \\
    Analytic Center & $\mathcal{O}({\cpDims} \log^2 ({\cpDims} \frac{R}{r}))$ &  $\mathcal{O}({\cpDims}^{1 + \omega} \log^2 ({\cpDims} \frac{R}{r}) + ({\cpDims} \log ({\cpDims} \frac{R}{r}))^{2 + \frac{\omega}{2}})$  & \cite{atkinson1995cutting}\\
    Random Walk & $\mathcal{O}({\cpDims} \log({\cpDims}\frac{R}{r}))$ & $\mathcal{O}({\cpDims}^7 \log({\cpDims}\frac{R}{r}))$ & \cite{bertsimas2004solving} \\
    Lee, Sidford, Wong & $\mathcal{O}({\cpDims} \log ({\cpDims} \frac{R}{r}))$ & $\mathcal{O}({\cpDims}^3 \log ^{\mathcal{O}(1)}({\cpDims} \frac{R}{r}))$ & \cite{lee2015faster}\\
    \hline
    \end{tabular}
    }
    \caption{This table summarizes the number of oracle queries and total runtime guarantees of a number of well-known cutting-plane variants. The stated runtimes are in terms of the problem dimension, $m$, the volume of the initial set, $R$, and the minimum volume of the set of interest, $r$. The constant $\omega$ represents the fast matrix multiplication constant.}
    \label{tab:cp}
\end{table}

\subsection{Application to Copositive Optimization}
Now that we have introduced cutting-plane algorithms, we are in a position to discuss their application to the copositive program \eqref{eq:COP}.
First, we will show how a \emph{copositivity oracle} can be readily transformed into a separation oracle for the feasible region of Problem \eqref{eq:COP}.
We will conclude with a discussion of how a copositivity oracle can be implemented using an Ising solver.
Formally, we define a copositivity oracle as follows:
\begin{definition}[Copositivity Oracle]
A copositivity oracle takes as input a matrix, $M$, and either returns \texttt{True} if $M$ is copositive or returns a vector $z \in \mathbb{R}^n_{\scriptscriptstyle \geq 0}$ such that $z^\top M z < 0$ (a ``certificate of non-copositivity").
\end{definition}

A copositivity oracle can be turned into a separation oracle for the feasible region of Problem \eqref{eq:COP} by expanding the terms in $z^\top M(\hat{\mu}, \hat{\lambda}, \hat{\gamma}) z$.
Explicitly, a test point, $(\hat{\mu}, \hat{\lambda}, \hat{\gamma})$, is infeasible if and only if $M(\hat{\mu}, \hat{\lambda}, \hat{\gamma})$ is not copositive.
Given $M(\hat{\mu}, \hat{\lambda}, \hat{\gamma})$ as input, the copositivity oracle returns a certificate of non-copositivity $z \in \mathbb{R}^{\mbqpNvars + 1}_{\scriptscriptstyle \geq 0}$ such that $z^\top M(\hat{\mu}, \hat{\lambda}, \hat{\gamma}) z < 0$.
In contrast, feasibility means that  $z^\top M(\mu, \lambda, \gamma) z \geq 0$.
Equivalently, the halfspace defined by
\begin{align}
    \shpConst &= z^\top \quadmat{Q}{c}{\cdot} z, \\
    \shpVec[\mu^{\text{(lin)}}_i] &= z^\top \quadmatrow{\cdot}{\frac{1}{2} A_{i, *}}{\cdot}z,\\
    \shpVec[\mu^{\text{(quad)}}_i] &= z^\top \Quadmat{A^\top_{i, *} A_{i, *}}{\cdot}{\cdot} z,\\
    \shpVec[\lambda_j] &= z^\top \quadmat{-\basis{j}  \basis{j} ^\top}{\frac{1}{2} \basis{j} }{\cdot} z,\\
    \shpVec[\gamma] &= z^\top \Quadmat{\cdot}{\cdot}{1} z,
\end{align}
is a separating hyperplane for $(\hat{\mu}, \hat{\lambda}, \hat{\gamma})$, where we use symbolic indexing to explicitly denote which variable each coefficient corresponds to.
Explicitly, the inner product between $a$ and $(\mu, \lambda, \gamma)$ is given by
\begin{equation}
    \shpVec^\top (\mu, \lambda, \gamma) = \sum_i \shpVec[\mu^{\text{(lin)}}_i]\mu^{\text{(lin)}}_i + \sum_i \shpVec[\mu^{\text{(quad)}}_i]\mu^{\text{(quad)}}_i + \sum_j \shpVec[\lambda_j]\lambda_j + \shpVec[\gamma] \gamma.
\end{equation}
This shows that given a copositivity oracle, constructing a separation oracle for Problem \eqref{eq:COP}, of dimension $\mathcal{O}(\mbqpNconstrs)$ and copositivity constraints on matrices of size $\mathcal{O}(\mbqpNvars)$, entails evaluating $\mathcal{O}(\mbqpNconstrs)$ vector-matrix-vector products, each of dimension $\mathcal{O}(\mbqpNvars)$.
The cutting-plane algorithms presented in Section \ref{subsec:cp_alg} can then be applied without further modification.

We note that the application of cutting-plane algorithms to copositive optimization has been explored from a classical perspective in \cite{guo2021copositive}, which considered their application to discrete markets and games, and \cite{badenbroek2022analytic}, which applied the algorithm to detect complete positivity of matrices. 
We believe our work is complementary to these prior works.
While our work relies on the off-the-shelf application of well-known cutting-plane variants, the algorithmic modifications in \cite{guo2021copositive}, and \cite{badenbroek2022analytic} provide insight for further improving our framework.
For example, the cutting plane algorithm in \cite{guo2021copositive} is readily hybridized leveraging our proposed approach.
We do not explore this extension in this paper due to their lack of convergence guarantees, while those of the variants presented are central to our theoretical analysis.
Moreover, the problem settings considered in these works serve as inspiration for additional applications that can be addressed with Ising solvers. 
On the other hand, the proof of strong duality (Theorem \ref{thm:strong_duality}) can be applied to address questions that were left open in \cite{guo2021copositive}.
We emphasize that the contribution of this work is not the copositive reformulation or a novel cutting-plane algorithm but rather the insight that these ideas are synergistic with recent advances in quantum(-inspired) computing.
In particular, copositive optimization is useful for deriving and analyzing new hybrid algorithms rooted in existing convex optimization algorithms, thus filling a gap in the hybrid algorithms literature.
\subsection{QUBO Approximation of Copositivity Checks}
Checking copositivity of $M(\mu, \lambda, \gamma)$ is naturally posed as the following (possibly non-convex) quadratic minimization problem
\begin{mini}
{z \in \mathbb{R}^{\mbqpNvars + 1}_{\scriptscriptstyle \geq 0}}{z^\top M(\mu, \lambda, \gamma) z}
{\label{eq:norm_cop}}{}
\addConstraint{||z||_p \leq 1,}
\end{mini}
where a matrix is copositive if and only if $\min \eqref{eq:norm_cop}$ is non-negative\footnote{While copositivity is defined as a condition over all of $\mathbb{R}^{\mbqpNvars + 1}_{\scriptscriptstyle \geq 0}$, quadratic scaling of the objective ensures that optimizing over a norm ball is sufficient for detecting copositivity.}.
There are several alternative approaches for checking copositivity \cite{anstreicher2021testing, dur2013testing, hiriart2010variational, bras2016copositivity, xia2020globally}; however, they are typically derived with Problem \eqref{eq:norm_cop} as the starting point and designed to exploit particular properties of Problem \eqref{eq:norm_cop}.
By choosing $p = \infty$, Problem \eqref{eq:norm_cop} can be approximated by a QUBO where an approximation of the matrix $M$, $\hat{M}$, is used such that the optimization variables $\hat{z}$ represent a binary expansion of $z$ with $k$ bits as follows:
\begin{mini}
{\hat{z}}{\hat{z}^\top \hat{M}(\mu, \lambda, \gamma) \hat{z}}
{\label{eq:qubo}\tag{QUBO}}{}
\addConstraint{\hat{z} \in \{0,1\}^{k(n+1)}}.
\end{mini}

Explicitly, $\hat{M}(\mu, \lambda, \gamma)$ and $M(\mu, \lambda, \gamma)$ are related as follows:
\begin{align}
    \hat{M}(\mu, \lambda, \gamma) &= \Dmat^\top M(\mu, \lambda, \gamma) \Dmat ,
\end{align}
where
\begin{align}
    \Dmat &:= \dfrac{1}{2^k - 1}
    \begin{pmatrix}
    2^0 & \cdots & 2^{\ndiscbin - 1} &
    0 & \cdots & 0 &
    \cdots & 
    0 & \cdots & 0\\
    0 & \cdots & 0 &
    2^0 & \cdots & 2^{\ndiscbin  - 1} &
    \cdots &
    0 & \cdots & 0\\
    \vdots & \vdots & \vdots &
    \vdots & \vdots & \vdots &
    \vdots &
    \vdots & \vdots & \vdots\\
    0 & \cdots & 0 & 
    0 & \cdots & 0 &
    \cdots &
    2^0 & \cdots & 2^{\ndiscbin  - 1} 
    \end{pmatrix},
\end{align}
The construction of \eqref{eq:qubo} is detailed in Appendix \ref{subsec:disc_cop}.
The explicit implementation of $\texttt{Oracle}(\cdot)$ is summarized in Algorithm \ref{alg:sep_or}.
Critically, the constraints of \eqref{eq:norm_cop} are implied by the natural domain of the Ising solver, mitigating the need to tune coefficients in a penalty method carefully.

\begin{algorithm2e}
\caption{Separation oracle, $\texttt{Oracle}(\cdot)$}\label{alg:sep_or}
\KwIn{$(\hat{\mu}, \hat{\lambda}, \hat{\gamma})$ (Test point)}
\KwOut{
\begin{equation*}
    \begin{cases}
    \texttt{True} & \text{ if } (\hat{\mu}, \hat{\lambda}, \hat{\gamma}) \text{ is feasible}\\
    \text{Separating hyperplane for } (\hat{\mu}, \hat{\lambda}, \hat{\gamma}) &  \text{ otherwise}
\end{cases}
\end{equation*}}
\tcp{Solve \eqref{eq:qubo} using an Ising solver}
\begin{argmini*}
{\hat{z}}{\eqref{eq:qubo}}
{}{z^*\gets}
\end{argmini*}

\uIf{$\min \eqref{eq:qubo} \geq 0$ }{
    \Return{\texttt{True}}\;
}
\Else{
    \begin{align}
        z &= \Dmat z^*\\
        \shpConst &= z^\top \quadmat{Q}{c}{\cdot} z \\
        \shpVec[\mu^{\text{(lin)}}_i] &= z^\top \quadmatrow{\cdot}{\frac{1}{2} A_{i, *}}{\cdot} z\\
        \shpVec[\mu^{\text{(quad)}}_i] &= z^\top \Quadmat{A^\top_{i, *} A_{i, *}}{\cdot}{\cdot} z\\
        \shpVec[\lambda_j] &= z^\top \quadmat{-\basis{j}  \basis{j} ^\top}{\frac{1}{2} \basis{j} }{\cdot} z\\
        \shpVec[\gamma] &= z^\top \Quadmat{\cdot}{\cdot}{1} z
    \end{align}
    \Return{$\shpVec, \, \shpConst$}\;
}
\end{algorithm2e}
\subsection{Discussion}
\begin{figure}
    \centering
    \includegraphics[width=\textwidth]{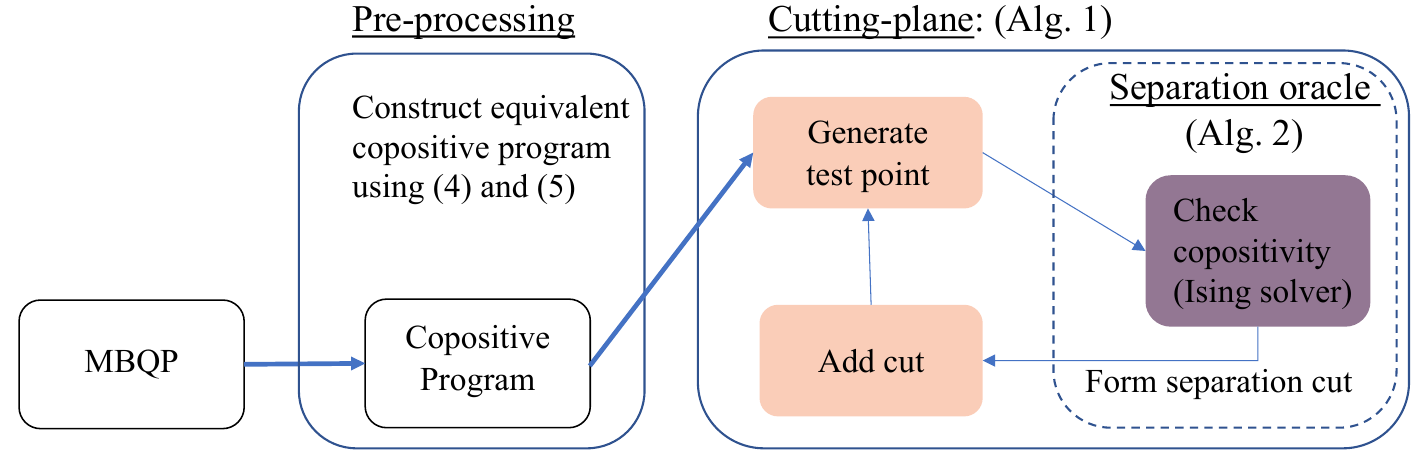}
    \caption{This figure depicts the entire solution process for solving a MBQP of the form \eqref{eq:MBQP}.}
    \label{fig:sol}
\end{figure}
In summary, we propose to solve Problem \eqref{eq:MBQP} by constructing the equivalent copositive formulation in \eqref{eq:COP} and applying any variant of Algorithm \ref{alg:cp_meta}. 
Within Algorithm \ref{alg:cp_meta}, the implementation of $\texttt{Oracle}(\cdot)$ is specified by Algorithm \ref{alg:sep_or}.
This process is depicted in Figure \ref{fig:sol}.
Now that we have presented our method in full, several comments are in order.
\paragraph{Computational complexity}
While the stated complexity of the cutting-plane algorithms is applicable to any problem, it is suggestively stated in terms of the variable $\cpDims$.
This notational overload is a deliberate choice because the dimension of the dual copositive program is equal to the total number of constraints in Problem \eqref{eq:CPP}, which is  $2\mbqpNconstrs + |B| + 1 = \mathcal{O}(\mbqpNconstrs)$.
The number of constraints can be reduced to $m + |B| + 1$ using the homogenized completely positive reformulation presented in Appendix \ref{subsec:strong_dual}.
While this will have no impact on the asymptotic complexity of the method, it can result in a practical reduction in runtime.
If $\mathcal{T}_Q$ represents the oracle complexity of a particular method, the total additional overhead of converting the copositivity oracle into a separation oracle is given by $\mathcal{O}(m n^2 \mathcal{T}_Q)$.
\paragraph{Discretization size}
Discretization of the copositivity check automatically introduces an approximation to the copositivity checks.
The approximation fidelity is improved as the number of discretization points is increased, although it is limited by the hardware.
Not only does representing a finer discretization require more qubits, but it also results in a greater skew in the coefficients of the Ising Hamiltonian.
This becomes challenging since many existing hardware platforms have limited precision in their implementable couplings.
In contrast, too coarse of a discretization runs the risk of missing the certificate of non-copositivity entirely.
This suggests that the discretization scheme should be well-tailored to the problem at hand; Appendix \ref{subsec:disc_cop} provides guidance for choosing a discretization size based on the coefficients of the Ising Hamiltonian.
A promising alternative is to circumvent discretization entirely, and apply quantum(-inspired) solvers that natively solve continuous variable box-constrained quadratic programs, such as the coherent continuous-variable machine (CCVM) recently proposed in \cite{khosravi2022non}.

\paragraph{Multiple cuts}
Following standard convention, this work assumes that the copositivity oracle returns a single value.
In contrast, in practice, many of the aforementioned Ising solvers are heuristics that involve multiple readouts.
Each of these reads can be used to construct a cut, where negative, zero, and positive Ising objective values correspond to deep, neutral, and shallow cuts, respectively.
Adding multiple cuts during each iteration is a possible heuristic for improving the convergence rate of the cutting-plane algorithm.
While the true ground state corresponds to the deepest cut, the convergence rate guarantees stated in Table~\ref{tab:cp} hold so long as a neutral or deep cut is added at each iteration.
Consequently, the proposed approach is not overly reliant on the Ising solver's ability to identify the ground state and is resilient to heuristics.
Critically, this raises the question of how to proceed if the Ising solver fails to return a certificate of non-copositivity, which will likely depend on problem specifics, such as the current outer approximation, the objective values of the samples, and the solver itself.
For example, if the Ising solver returns positive but small solutions, depending on the current outer approximation, the addition of shallow cuts can still reduce the search space.
On the other hand, if all non-zero solutions result in a large objective value, one could increase confidence that the test point is feasible by increasing the number of discretization points and readouts.

\section{Experiments}\label{sec:experiments}

We conducted an investigation of the proposed method on the maximum clique problem, which finds the largest complete subgraph of a graph.
Given a graph, the maximum clique problem can be formulated as a completely positive program 
\begin{maxi}
    {X \in \mathbb{R}^{n \times n}}{\innerMat{\mathds{1} \mathds{1}^\top}{X}}
    {}{\label{eq:mc_cpp}}
    \addConstraint{\innerMat{\overline{A} + I}{X} = 1}
    \addConstraint{X \in \mathcal{C}^*_n,}
\end{maxi}
where $\overline{A}$ is the adjacency matrix of the graph's complement \cite{de2002approximation}.
We note that solving the maximum clique problem is equivalent to solving the maximum independent set problem on the complement graph.
The dual of \eqref{eq:mc_cpp} is the following copositive program:
\begin{mini}
    {\lambda \in \mathbb{R}}{\lambda}
    {}{\label{eq:mc_cop}}
    \addConstraint{\lambda(I + \overline{A}) - \mathds{1} \mathds{1}^\top \in \mathcal{C}_n.}
\end{mini}
This copositive program only has one variable regardless of the graph's number of vertices or edges.
Thus, we solve it with bisection, a special case of the ellipsoid algorithm, as the cutting-plane algorithm.
The copositivity check's size is determined by the number of vertices, $n$, which impacts the complexity of computing the cuts from the certificates of non-copositivity.
The number of edges can be used to upper-bound the size of the maximum clique, thus determining the size of the initial feasible region; however, its effect on the complexity of checking copositivity is unclear.

\begin{figure}
     \centering
     \begin{minipage}[b]{.2\textwidth}
     \subfloat[]
    {\label{fig:mc_ex}\includegraphics[clip, trim=6cm 8cm 5cm 3cm,width=\linewidth]{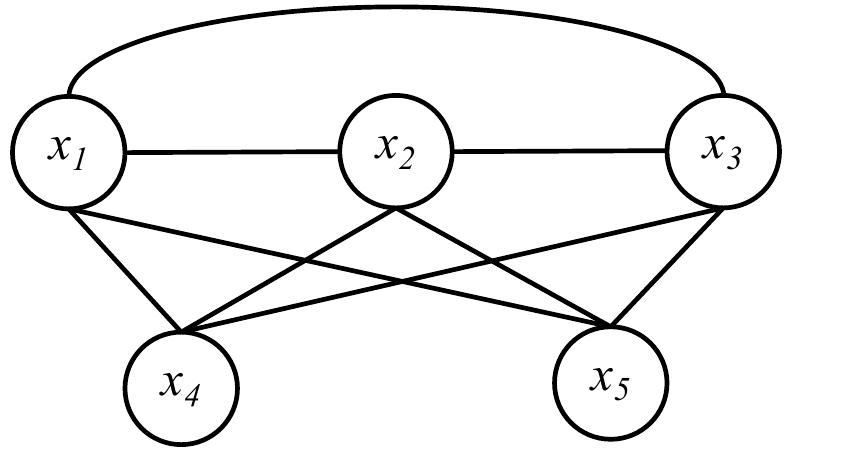}
    }
    \vfill
     \subfloat[]
     {\label{fig:mc_ex_mat}\includegraphics[width=\linewidth]{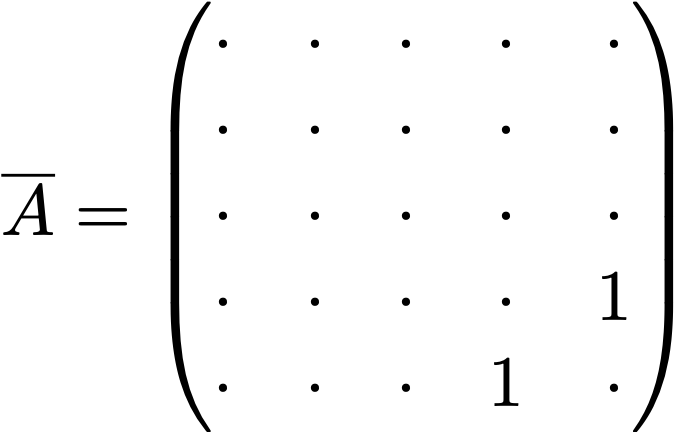}}
     \end{minipage}
     \hfill
    \begin{minipage}[b]{0.75\textwidth}
     \subfloat[]
        {\label{fig:mc_ex_sol}\includegraphics[width=\linewidth]{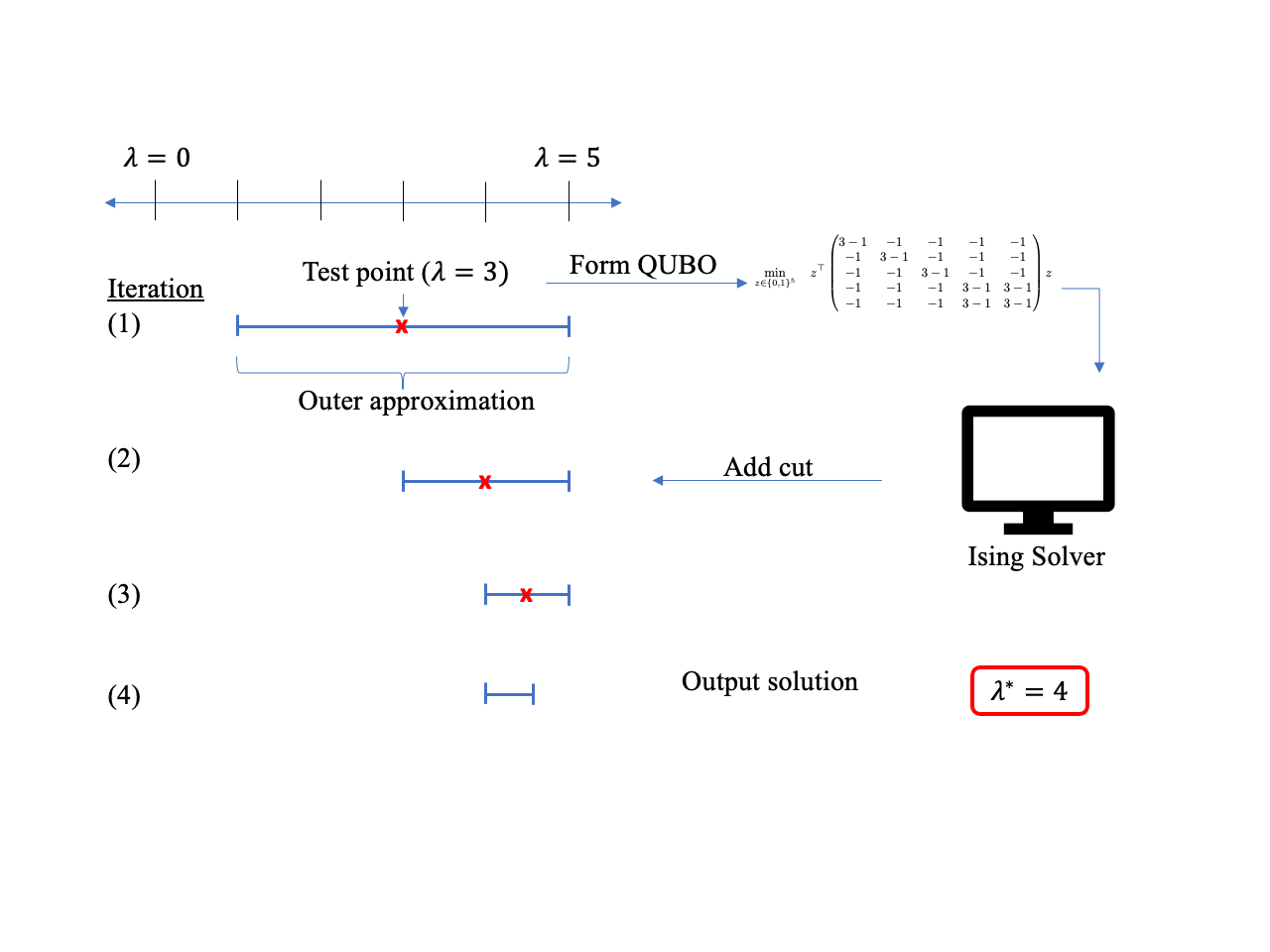}}
     \end{minipage}
     \caption{Figure \protect\ref{fig:mc_ex} depicts a small maximum clique example where there are edges between all vertices except $x_4$ and $x_5$. Figure \protect\ref{fig:mc_ex_mat} depicts the adjacency matrix of graph \protect\ref{fig:mc_ex}'s complement, which has a single edge between vertices $x_4$ and $x_5$. Figure \ref{fig:mc_ex_sol} depicts the solution process for the copositive cutting-plane algorithm.}
\end{figure}

\subsection{QUBO Subroutine}
\begin{figure}
    \centering
    \includegraphics[width=\textwidth]{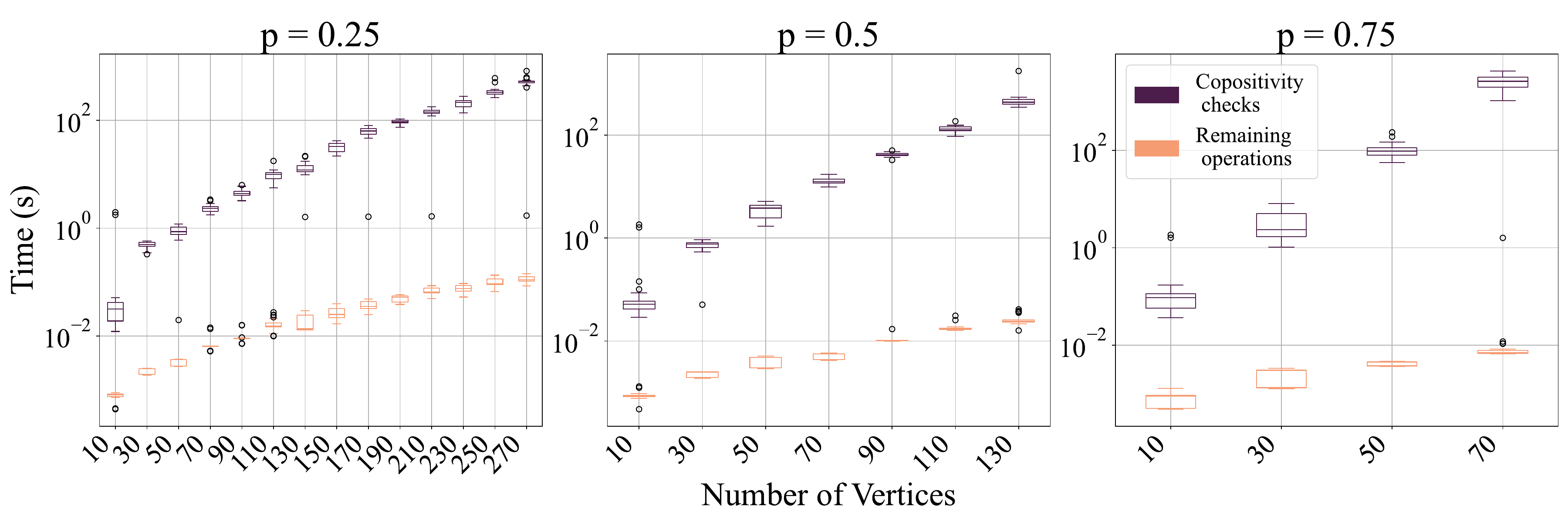}
    \caption{This figure plots the time spent on the copositivity checks versus all other operations in the proposed method. The copositivity checks grow exponentially with the number of vertices, while the other operations grow modestly.}
    \label{fig:profile}
\end{figure}

Many hybrid algorithms are designed by replacing subroutines of existing (fully classical) algorithms with a quantum(-inspired) counterpart.
An oft-neglected consideration is whether the quantum(-inspired) computer is applied to a bottleneck in the original algorithm.
We contend that for a hybrid algorithm to yield significant speed-ups, the quantized subroutine should constitute the bulk of the algorithm's complexity. 
In keeping with this supposition, we first evaluated whether the copositive cutting-plane algorithm shifts the complexity of the solution process onto the copositivity checks by profiling each component of the algorithm separately.

To study the scaling of each component of the proposed approach, we considered random max-clique problems with up to $10, 30, \ldots, 270$ vertices, where the maximum graph size varies by edge density to ensure a reasonable computation time for this experiment.
For each graph size, we generated 25 random Erd\H{o}s-Renyi instances with edge densities $p \in \{0.25, \, 0.5, \, 0.75 \}$ and solved to an absolute gap of 0.9999 between upper and lower bounds (because the optimal solution is known to be integral) using the proposed copositive cutting plane algorithm.
The copositivity checks were conducted by solving Anstreicher's MILP characterization of copositivity~\cite{anstreicher2021testing} (which we found to be one of the most competitive classical formulations), using \texttt{Gurobi} version 9.0.3 \cite{GurobiOptimization2020}.
\footnote{All experiments were run on an AMD Ryzen 7 1800X Eight-Core Processor@3.6GHz with 64GB of RAM and 16 threads.
All code needed to reproduce these experiments is available at \href{https://github.com/StanfordASL/copositive-cutting-plane-max-clique}{https://github.com/StanfordASL/copositive-cutting-plane-max-clique}.}

Figure~\ref{fig:profile} plots the time the copositive cutting plane algorithm spent on the copositivity checks versus other operations (updating the outer approximation and computing test points).
The time spent on the copositivity checks scales exponentially with the number of vertices in the graph, while the time spent on other operations grows modestly.
This is because Problem \eqref{eq:mc_cpp} only has one constraint regardless of the graph's number of vertices or edges. 
In contrast, the size of the copositivity check is exactly equal to the number of vertices in the graph.
Both the theoretical analysis and empirical results confirm that the proposed approach shifts the complexity of the copositive program onto the copositivity checks.
This experiment shows that the proposed methodology is particularly effective for problems whose constraints remain constant or grow modestly with problem size.

\begin{figure}[t!]
    \centerline{\includegraphics[width=0.95\textwidth]{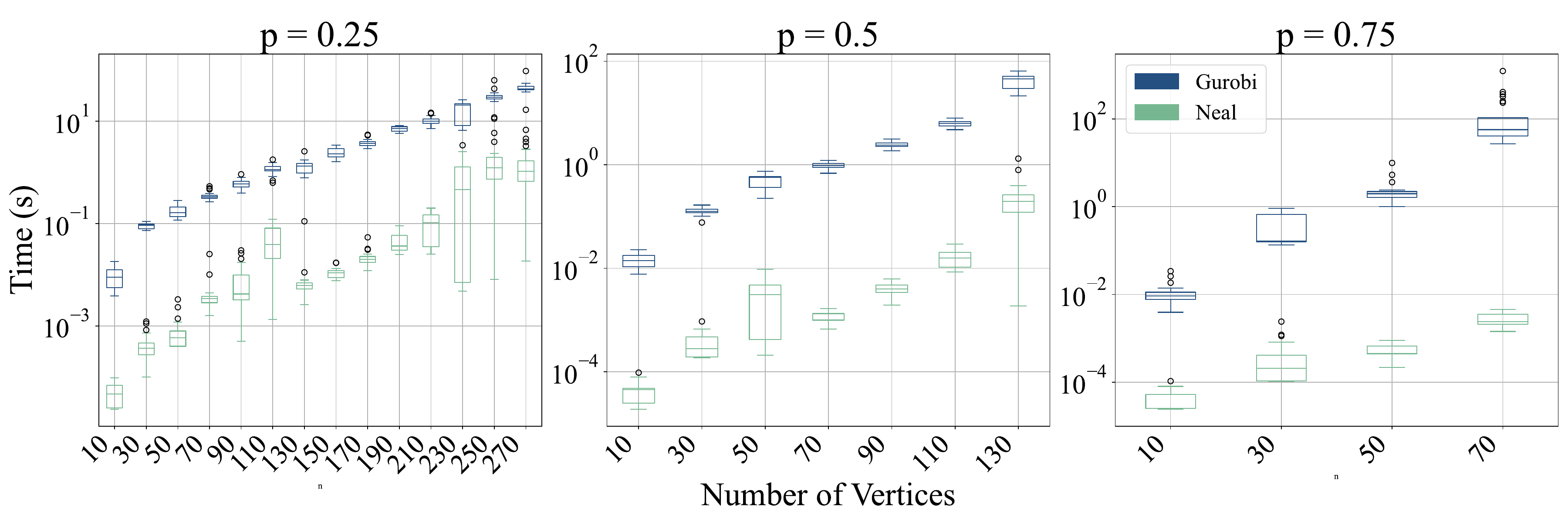}}
    \caption{This figure plots the time to target to 99\% confidence of the Simulated Annealing (SA) implementation in \texttt{Neal} (replacing the continuous feasible region, $[0, 1]^n$, with its vertices, $\{0, 1\}^n$) against the solution time of \texttt{Gurobi} for the copositivity checks.
    We solved each copositivity check with 100 sweeps and 1000 reads.
    For all densities, both methods scale exponentially with the number of vertices in the graph; however, SA is several orders of magnitude faster than \texttt{Gurobi}.}
    \label{fig:DG_comp}
\end{figure}

While undesirable for a fully classical implementation, the overwhelming complexity in the copositivity checks represents an opportunity for a quantum(-inspired) solver to beget significant speed-up.
They will result from being able to execute the copositivity checks faster than the classical implementation (i.e., Anstreicher's MILP formulation).
To investigate potential speedups from using a stochastic Ising solver, we re-solved each copositivity check that yielded a certificate of non-copositivity using Simulated Annealing (SA) through the software \texttt{Neal} version 0.5.9, a SA sampler 
\cite{DWaveNeal2021}, i.e., a solver that returns samples on the solutions distribution generated by SA.
Because SA is not guaranteed to find the global optima in a single annealing cycle, we define a probabilistic notion of time to target.
In particular, we follow \cite{ronnow2014defining} and define the time to target with $s$ confidence to be the number of repetitions to find the ground state at least once with probability $s$ multiplied by the time for each annealing cycle, $\texttt{T}_{\text{anneal}}$, i.e., 
\begin{align}
    \texttt{TTT}_s = \texttt{T}_{\text{anneal}}\frac{\log(1 - s)}{\log(1 - \hat{p}_{\text{succ}})}, 
\end{align}
where $\hat{p}_{\text{succ}}$ is the expected value of the returned solution divided by the ground state/minimum.
This results in a probability of success that interpolates between counting only solutions corresponding to the ground state and counting all certificates of non-copositivity as successes by considering the relative quality of each sample.
We will also consider analogous scenarios where only ground state solutions are counted as success; we reserve the terminology ``time to solution", $\texttt{TTS}_s = \texttt{T}_{\text{anneal}}\frac{\log(1 - s)}{\log(1 - p_{\text{succ}})}$, for such cases to distinguish from the previously defined time to target. 
The values of $\hat{p}_{\text{succ}}$ and $p_{\text{succ}}$ is evaluated empirically over 1000 samples/reads.
The time per annealing cycle, $\texttt{T}_{\text{anneal}}$, was evaluated as the total wall-clock time (for all reads) divided by the number of reads.
All other \texttt{Neal} parameters were left as their default values.

The time to solution and time to target metrics are intended to facilitate comparison between deterministic and stochastic solvers by taking into account both the time needed to run the stochastic solver and its probability of success.
It is impossible to guarantee 100\% success for a stochastic solver under this formulation (mathematically, this would be equivalent to trying to compute $\texttt{TTS}_1$ or $\texttt{TTT}_1$, which are both undefined). 
However, computing these metrics with a high degree of certainty (e.g., $s = 1 - \epsilon$) is widely accepted as a tolerable, albeit imperfect, benchmark \cite{king2015benchmarking}.
In the remainder of this section, we will compare time to target metrics from SA against solution time from \texttt{Gurobi}, but the astute reader should keep in mind that the two solvers serve fundamentally different purposes.
This comparison between the solvers is not intended to be interpreted in isolation, but rather to highlight where it might be appropriate and beneficial to substitute a heuristic Ising solver within the copositive cutting-plane framework.

For each copositivity check solved, we considered discretizations corresponding to \linebreak ${\min_{\hat{z} \in \{0, \, 1\}^n} \hat{z}^\top M \hat{z} }$ (i.e., no additional problem discretization).
We solved each copositivity check with 100 sweeps and 1000 reads\footnote{While the performance of \texttt{Neal} depends on the number of sweeps, we found that optimizing the number of sweeps does not result in significant reductions in the time to target. 
We provide further discussion in the Appendix \ref{subsec:hyperopt}.}.
Figure \ref{fig:DG_comp} plots the time to target with 99\% confidence from \texttt{Neal} against the solution time from \texttt{Gurobi}\footnote{Note that \texttt{Gurobi}'s solution time in Figure \ref{fig:DG_comp} is different from copositivity checks profiling in Figure \ref{fig:profile}. This is because only non-copositive instances were considered for this comparison, while all instances, including copositive ones, were included in the profiling comparison.}.
We see that for all graph sizes, \texttt{Neal} can consistently find certificates of non-copositivity in orders of magnitude less time than \texttt{Gurobi}.
Notably, \texttt{Neal} and \texttt{Gurobi} demonstrate similar scaling with respect to the number of vertices.

\begin{figure}
    \centering
    \includegraphics[width=\textwidth]{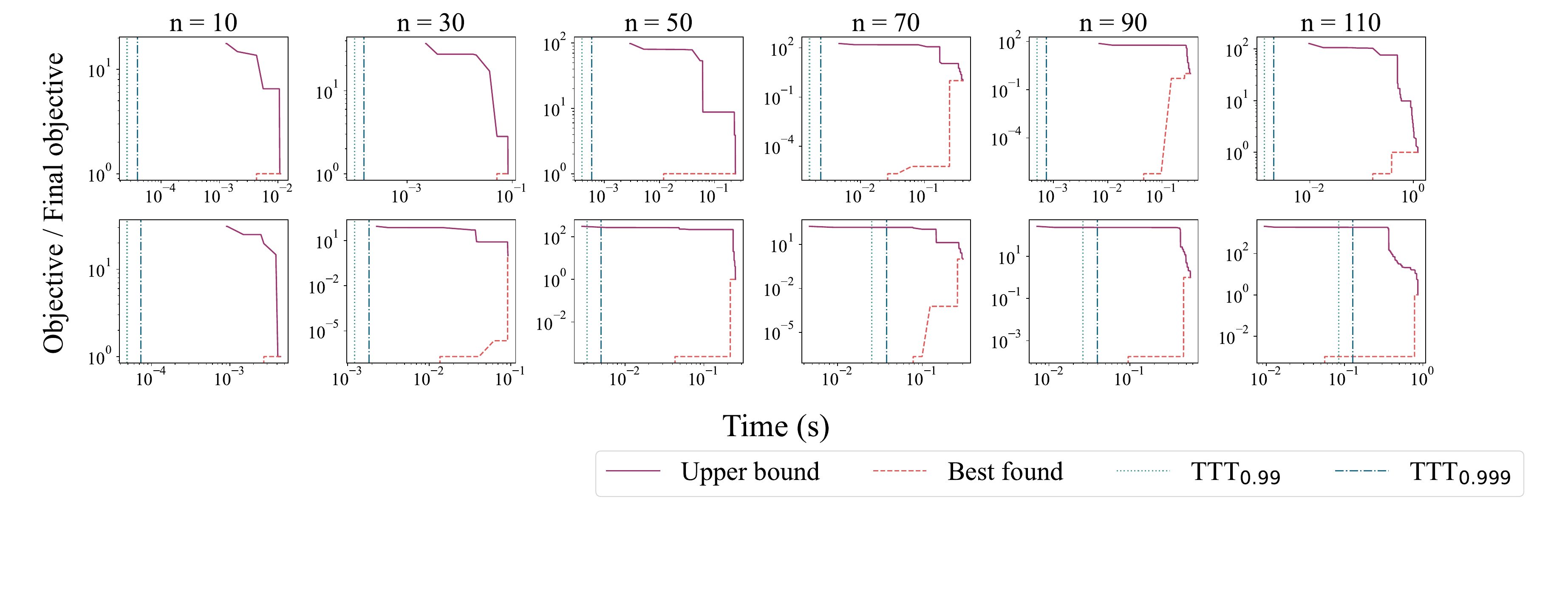}\par\vspace{-2\baselineskip}
    \caption{This figure depicts sample trajectories of \texttt{Gurobi}'s upper and lower bounds against $\texttt{TTT}_{0.99}$ and $\texttt{TTT}_{0.999}$ for edge density $p=0.25$.
    For each graph size, the top row represents the instance where the ratio between \texttt{Gurobi}'s solution time and $\texttt{TTT}_{0.99}$ is the greatest, and the bottom row represents the instance where the ratio is the smallest--all instances were run with 100 sweeps.
    In most instances, \texttt{Neal} reaches the $\texttt{TTT}_{0.999}$ confidence before \texttt{Gurobi} even returns a callback.}
    \label{fig:min_max_gurobi}
\end{figure}

Unlike SA, which operates without reference to rigorous optimality bounds, \texttt{Gurobi}'s solution process tracks both upper and lower bounds on the objective value and terminates only when they reach user-specified stopping conditions.
To evaluate whether the optimal objective is found early in the solution process and time is spent closing the upper bounds, we plotted \texttt{Gurobi}'s lower and upper bounds progress against time together with $\texttt{TTT}_{0.99}$ and $\texttt{TTT}_{0.999}$ in Figure \ref{fig:min_max_gurobi} for instances with density $p = 0.25$.
Analogous plots for other densities are included in the Appendix.
For each graph size, we plotted the instances where the ratio between \texttt{Gurobi}'s solution time and $\texttt{TTT}_{0.99}$ is the greatest (top row) and least (bottom row)--all instances were run with 100 sweeps.
For each instance, we plot \texttt{Gurobi}'s upper bound (purple, solid) and best objective found (red, dashed), and \texttt{Neal} $\texttt{TTT}_{0.99}$ (light teal, dashed), and  $\texttt{TTT}_{0.999}$ (dark teal, dash-dot) (colored figures are available online).
We found that in most instances, \texttt{Neal} reaches the time to target with 99.9\% confidence before \texttt{Gurobi} even returns a callback (i.e., when the dark teal, dash-dot line does not intersect either of the purple-solid or red-dashed lines); this is likely due to an initial pre-processing step.
Pre-processing is a necessary overhead for Gurobi's intended purpose of proving optimality, thus hampering its relative performance on smaller or easier problems.
Critically, in the proposed approach, optimality guarantees are only necessary for cases where the test point is copositive, while the algorithm can make progress with any certificate of non-copositivity, even if it is not globally optimal.
This suggests that heuristics (e.g., Ising solvers) and complete methods (e.g., \texttt{Gurobi}) could play complementary roles within the same copositive cutting-plane algorithm.
For example, for some test points, one may avoid the overhead of \texttt{Gurobi} altogether if the Ising solver quickly returns a certificate of non-copositivity.
On the other hand, if the Ising solver fails to generate such a certificate, one may rely on \texttt{Gurobi} for proving that the test point is copositive.

\subsection{Overall Cutting-Plane Algorithm}
\begin{figure}
    \centering
    \includegraphics[width=\textwidth]{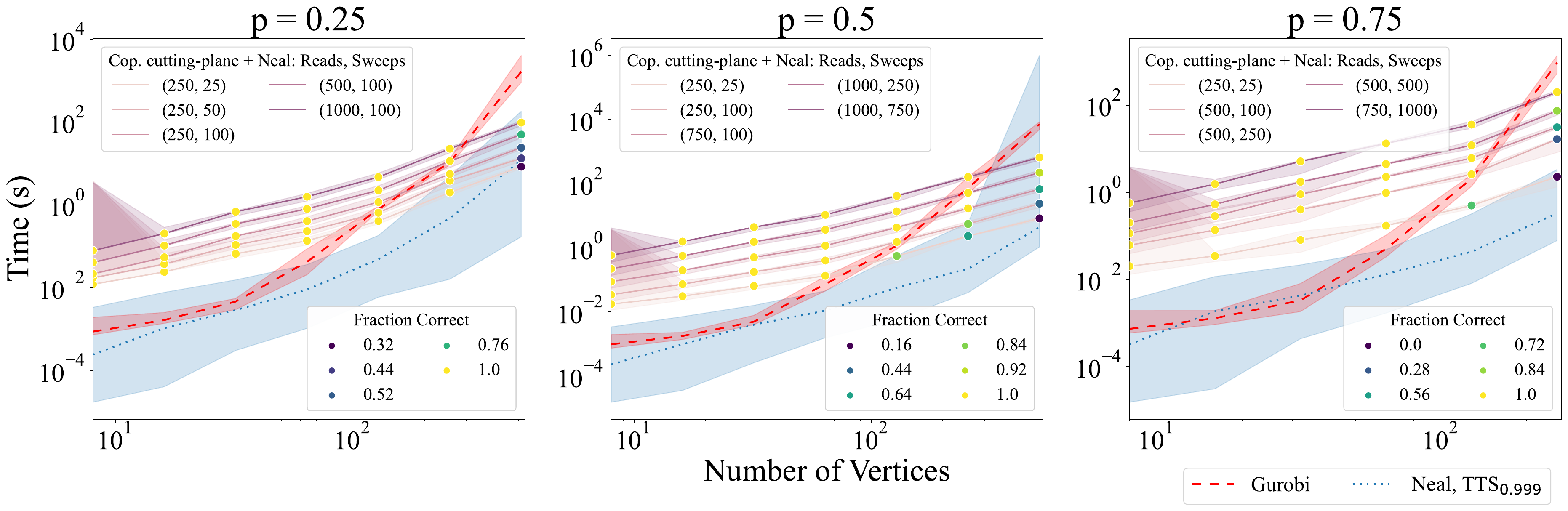}
    \caption{This figure plots the solution time for the copositive cutting-plane algorithm with the Simulated Annealing implementation in \texttt{Neal} as the Ising solver (with various fixed parameters), the solution time when solving a mixed-integer programming (MILP) formulation of maximum-clique directly with \texttt{Gurobi}, and the corresponding $\texttt{TTT}_{0.999}$ from \texttt{Neal} applied to the \texttt{maximum\_clique} formulation with penalty weight 1.
    While \texttt{Gurobi}'s solution time is orders of magnitude faster than the copositive cutting-plane algorithm for the smallest graph sizes, the copositive cutting-plane algorithm starts outperforming it for larger graph sizes.
    }
    \label{fig:all_comparison}
\end{figure}

Next, we compared the copositive cutting-plane algorithm with the SA implementation in \texttt{Neal} as the Ising solver against solving a mixed-integer linear program (MILP) formulation of maximum-clique directly with \texttt{Gurobi}. 
\texttt{Gurobi}'s solution time was evaluated on the following MILP formulation of maximum clique:
\begin{maxi}
{x \in \{0, 1\}^n}{\mathds{1}^\top x}
{}{}
\addConstraint{x_i + x_j \leq 1, \quad}{\forall (i, j) \in \overline{E},}
\end{maxi}
where $\overline{E}$ is the edges in the complement graph.
This is a MILP with $n$ binary variables, where $n$ is the number of vertices in the graph, and $|\overline{E}|$ constraints (i.e., the number of edges in the complement graph).
The copositive cutting-plane algorithm was tested with different sweeps and reads, which were fixed throughout each run of the algorithm.
The solid pink lines in Figure \ref{fig:all_comparison} plot the runtime of the copositive cutting-plane algorithm for a representative set of these parameters.
Because \texttt{Neal} may fail to find a certificate for some non-copositive matrices, this method may incorrectly reduce the upper bound in the outer approximation; however, it cannot incorrectly update the lower bound.
Throughout the algorithm, we track the exact lower bound and the approximate upper bound (which is updated when \texttt{Neal} fails to find a certificate of non-copositivity).
The algorithm is terminated when the lower bound and approximate upper bound are within an absolute tolerance of 0.9999; this is the same stopping condition as the exact case where we checked copositivity using Anstreicher's MILP formulation.
Consequently, the solution returned was determined by rounding the lower bound up to the nearest integer.
This means that the solution returned is guaranteed to be a lower bound for the maximum clique instance.
The fraction of correct maximum clique solutions is indicated by the color of the markers.

From the pink plots, we see that for a \emph{fixed parameter setting}, the copositive cutting-plane algorithm exhibits polynomial scaling in the graph size.
While it is tempting to extrapolate this scaling relationship to larger graph sizes, the failure of parameters that were successful on smaller graphs on larger graph sizes (denoted by the green and blue markers) indicate that it is unlikely that fixed parameter settings will continue to be effective \emph{ad infinitum}.
In particular, we found that while smaller and sparser instances are tolerant of fewer sweeps and reads (resulting in shorter calls to the Ising solver), larger instances required more sweeps and reads to be accurate. 
While the copositive cutting-plane algorithm is designed to benefit from speed-ups of state-of-the-art Ising solvers, the converse is also expected; poor scaling of the Ising solver will make its way into the runtime as the time spent in each oracle call, $\mathcal{T}_Q$.

On the other hand, we observe that the confidence intervals for \texttt{Neal} are significantly wider than those of \texttt{Gurobi} and the copositive cutting-plane algorithm.
This is potentially due to the sensitivity of \texttt{Neal} to its penalty weight.
While \texttt{Neal} can also produce lower bounds for the maximum clique, the bounds can only be updated when it returns solution corresponds to a clique.
In contrast, the copositive cutting-plane algorithm can generate cuts from any certificate of non-copositivity, even those that do not correspond to a clique at all.
This means that the copositive cutting-plane algorithm can make progress from a larger set of the solutions returned by the Ising solver. 
This suggests that the copositive cutting-plane algorithm may even be competitive against its underlying copositivity checker solving a direct formulation of the problem, particularly if its performance is highly sensitive to parameter settings.

We note that \texttt{Gurobi} takes advantage of multi-threading while neither \texttt{Neal} nor the copositive cutting-plane algorithm do.
This raises the important question of how much the copositive cutting-plane algorithm (and \texttt{Neal}) will benefit from similar decomposition and parallelization efforts.

\begin{sidewaysfigure}
    \centering
    \includegraphics[width=\textheight]{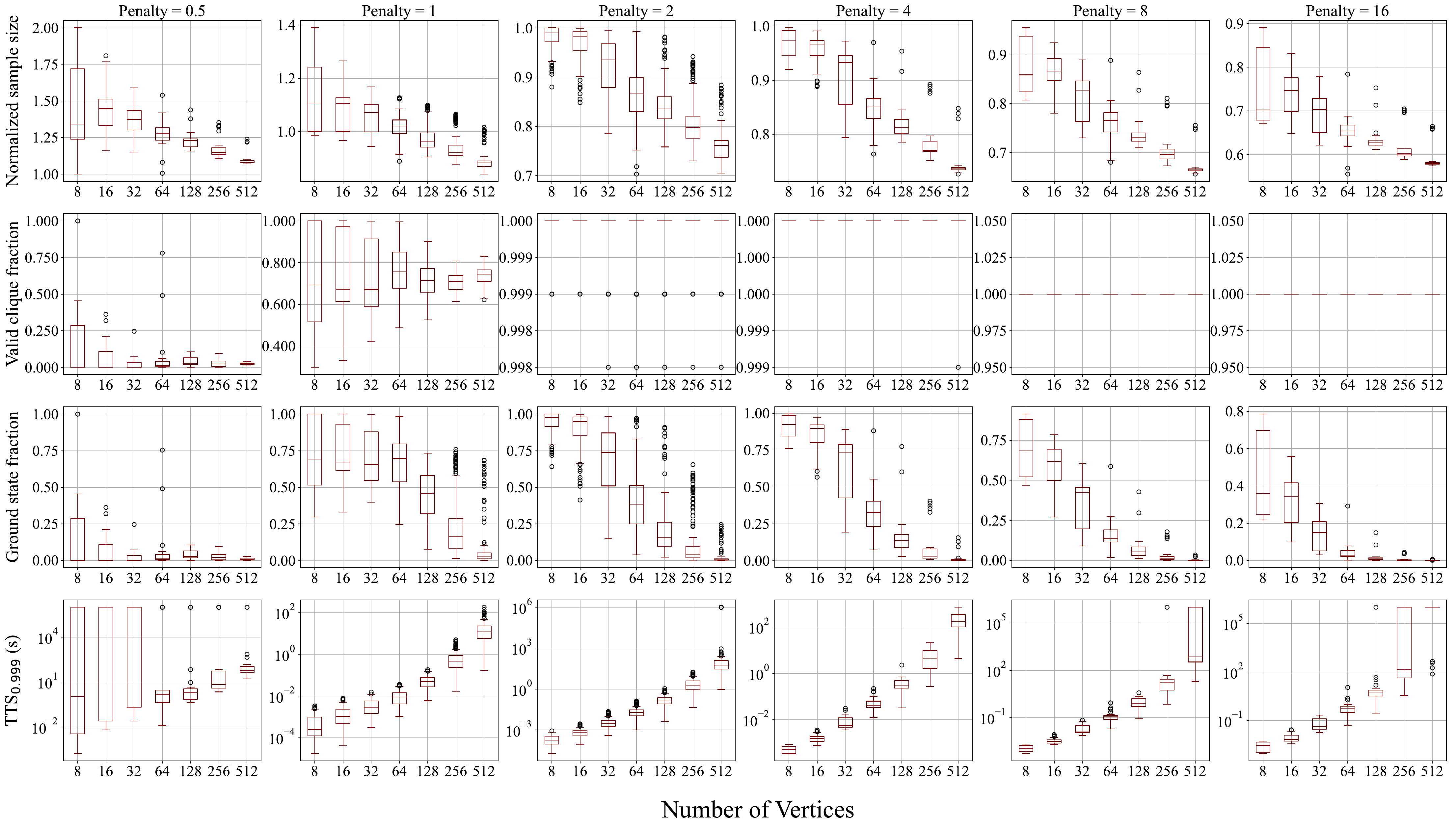}
    \caption{This figure plots the normalized sample size (the size of the returned solution divided by the ground truth maximum clique size) and the fraction of reads that resulted in a valid clique for graph density $p = 0.25$. These figures were used to compute the fraction of reads resulting in a ground state solution and the corresponding $\texttt{TTT}_{0.999}$ (also plotted). As the penalty weight is increased, the normalized sample size decreases, and the fraction of valid cliques increases. This highlights the delicate trade-off between constraints and the objective in penalty formulations.}
    \label{fig:penalty_25}
\end{sidewaysfigure}

Finally, we investigated the effectiveness of directly converting the maximum clique problem to an Ising problem using a standard penalty formulation.
To do so, we solved each of the maximum clique problem instances using the \texttt{maximum\_clique} formulator\footnote{\url{https://docs.ocean.dwavesys.com/projects/dwave-networkx/en/latest/reference/algorithms/generated/dwave_networkx.maximum_clique.html}} with \texttt{Neal} as the sampler and a range of penalty weights in $\{2^{-1}, 2^{0}, \ldots, 2^4\}$; the number of sweeps was left to its default value of 1000.
This results in a QUBO with $n$ variables and $|\overline{E}|$ quadratic terms.
For each instance, we conducted 1000 reads and evaluated the average normalized sample size (the size of the returned solution divided by the ground truth maximum clique size) and the fraction of reads that resulted in a valid clique; a ground state solution is one that is both a valid clique and has a normalized sample size of 1.
We computed the probability of success, $p_{\text{succ}}$, as the fraction of reads that resulted in a ground state solution, which was subsequently used to derive the time to solution to 99.9\% confidence.
Figure \ref{fig:penalty_25} plots each of these metrics as a function of the penalty weights and graph size for edge density $p = 0.25$.
Analogous plots for other densities are included in the Appendix.

For penalty weights $0.5$ and $1$, the normalized sample size is often greater than $1$, resulting in samples that do not represent a valid clique.
For penalty weights $2, 4, 8$, and $16$, most samples were valid cliques; however, the normalized sample sizes were typically less than 1; these represent non-maximum cliques.
Generally, as the penalty weight is increased, the normalized sample size decreases, and the fraction of valid cliques increases.
This aligns with the interpretation that the penalty weight represents a trade-off between satisfying the constraints versus optimizing the objective.
These empirical results also corroborate the analytical results of \cite{quintero2022characterization}, which state that the minimum valid penalty weight for the stable set of a graph is 1.
Given that \texttt{maximum\_clique} represents the maximum clique problem as finding the stable set of the graph built with the complement of the original edges, the bound on the penalty weight is valid.
This experiment demonstrates that while the penalty formulation may be an effective heuristic, it typically requires carefully tuning the penalty weights to optimize the trade-off between satisfying the constraints and optimizing the objective.
Figure \ref{fig:all_comparison} also plots the corresponding $\texttt{TTT}_{0.999}$ from \texttt{Neal} applied directly to the \texttt{maximum\_clique} formulation with penalty weight 1 against the \texttt{Gurobi} solution time and the copositive cutting-plane solution time.

\section{Conclusions}\label{sec:conclusion}
In this paper, we advocate for the development of a theory of hybrid quantum-classical algorithms that analytically quantifies their performance.
Metrics for comparing different hybrid algorithms may include the number of calls to the quantum computer, the complexity of the classical portion, and the requirements on the quantum computer.
As a step in this direction, we demonstrate a class of hybrid algorithms for mixed-binary quadratic programming problems using Ising solvers and report the aforementioned metrics.
Our framework relies on Burer's convex reformulation of such problems using completely positive programming.
Our first contribution is to extend this result and show that under mild conditions, the dual copositive program exhibits strong duality.
We then propose a hybrid quantum-classic solution algorithm based on cutting-plane algorithms, where an Ising solver is used to construct the separation oracle.
This approach partially mitigates the heuristic nature of many state-of-the-art Ising solvers.
Moreover, the runtime of the components handled by the classical computer scales polynomially with the number of constraints in the original mixed-binary quadratic program.
This suggests that if our approach is applied to a problem with exponential scaling, the complexity is shifted on the subroutine carried out by the hardware accelerator, e.g., the quantum computer.
Our proposed approach is particularly appealing because it suggests that the proposed approach could take advantage of any speedup that exists even without an explicit characterization of what that speedup is.

While the proposed framework seems like a promising way forward for utilizing quantum/quantum-inspired Ising solvers, a crucial question remains regarding how the algorithm should proceed if the Ising solver fails to find a certificate of non-copositivity.
Could one design an efficient algorithm that circumvents the ambiguity due to failure to find a certificate of non-copositivity (perhaps using a sum-of-squares-based inner approximations of the copositive cone)?
Alternatively, is there a complexity barrier that prevents such a construction?
More broadly, identifying fundamental limitations, such as this example, is important for understanding what requirements should be placed on new computing architectures. 
While the theory in this paper is framed in the context of understanding algorithms that interact with existing hardware, we believe its most potent impact is in informing the co-design of ``hardware primitives" and optimization algorithms.
As a concrete example, the coherent continuous variable machine (CCVM) recently proposed in \cite{khosravi2022non} circumvents discretization of the copositivity checks entirely, thus reducing the hardware resources and potentially improving problem conditioning for each copositivity check.
Not only does this hardware innovation directly benefit the optimization community, but the analysis in this work helps to justify the development of continuous variable devices.

\subsection*{Acknowledgements}
This work was supported by NSF CCF (grant \#1918549), NASA Academic Mission Services (contract NNA16BD14C – funded under SAA2-403506).  
R.B. acknowledges support from the NASA/USRA Feynman Quantum Academy Internship program.
The authors wish to thank Aaron Lott, Luis Zuluaga, and Juan Vera for helpful input and discussions during the development of these ideas.

\bibliographystyle{IEEEtran}
\bibliography{references.bib}

\section{Appendix}
\subsection{Proof of strong duality}\label{subsec:strong_dual}
Problem \eqref{eq:CPP} is equivalent to the following homogenous form completely positive program (i.e., $ \min \eqref{eq:CPP} = \min \eqref{eq:CPP_hom}$):
\begin{mini}
    {x \in \mathbb{R}^{\mbqpNvars},X \in \mathbb{R}^{\mbqpNvars \times \mbqpNvars}}{\innerMat{\quadmat{Q}{c}{\cdot}}{\quadmat{X}{x}{1}}}
    {\label{eq:CPP_hom}\tag{Hom-CPP}}{}
    \addConstraint{\innerMat{\quadmatrow{A^\top_{i, *} A_{i, *}}{-b_i A_{i, *}}{b_i^2}}{\quadmat{X}{x}{1}} = 0}
    \addConstraint{\innerMat{\quadmat{-\basis{j}  \basis{j} ^\top}{\frac{1}{2} \basis{j} }{\cdot}}{\quadmat{X}{x}{1}} = 0, \, \forall j \in B}
    \addConstraint{\quadmat{X}{x}{1} \in \mathcal{C}^*_{\mbqpNvars + 1}.}
\end{mini}
This form will be useful for proving strong duality.
Because the homogenized form of the equality constraints form a cone, this perspective will help prove strong duality between Problem \eqref{eq:CPP} and its dual.
The Lagrangian dual of \eqref{eq:CPP_hom} is the following copositive optimization problem:
\begin{maxi}
{\mu,\lambda,\gamma}{\gamma}
{\label{eq:COP_hom}\tag{Hom-COP}}{}
\addConstraint{\hat{M}(\mu, \lambda, \gamma) \in \mathcal{C}_{\mbqpNvars + 1},}
\end{maxi}
where $\hat{M}(\mu, \lambda, \gamma)$ is defined as 
\begin{equation}
    \begin{aligned}
    \hat{M}(\mu, \lambda, \gamma) := &\quadmat{Q}{c}{\cdot} - \sum_i \mu_i \quadmatrow{A^\top_{i, *} A_{i, *}}{-b_i A_{i, *}}{b_i^2}\\
    &- \sum_{j \in B} \lambda_j \quadmat{-\basis{j} \basis{j}^\top}{\frac{1}{2} \basis{j}}{\cdot} - \gamma\Quadmat{\cdot}{\cdot}{1}.
    \end{aligned}
\end{equation}

\begin{theorem}[Homogeneous Strong Duality]\label{thm:hom_strong_duality} If Problem \eqref{eq:MBQP} is feasible with bounded feasible region, then strong duality holds between Problems \eqref{eq:CPP_hom} and \eqref{eq:COP_hom}.
Moreover, an $\epsilon$ optimal value of \eqref{eq:COP_hom} is obtained by a feasible solution, where $\epsilon > 0$ can be arbitrarily small.
\begin{proof}
Notice that the set of affine constraints, 
\begin{equation}
    \mathcal{NULL} := \left\{ \tilde{X} \mid \innerMat{\quadmatrow{A^\top_{i, *} A_{i, *}}{-b_i A_{i, *}}{b_i^2}}{\tilde{X}} = 0, \, \innerMat{\quadmat{-\basis{j} \basis{j}^\top}{\frac{1}{2} \basis{j}}{\cdot}}{\tilde{X}} = 0\right\},
\end{equation}
forms a cone.
So, we could express \eqref{eq:CPP_hom} as the following optimization problem:
\begin{mini}
    {\tilde{X} \in \mathbb{R}^{(\mbqpNvars + 1) \times (\mbqpNvars + 1)}}{\innerMat{\quadmat{Q}{c}{\cdot}}{\tilde{X}}}
    {}{\label{eq:CPP_2cone}}
    \addConstraint{\innerMat{\Quadmat{\cdot}{\cdot}{1}}{\tilde{X}} = 1}
    \addConstraint{\tilde{X} \in \mathcal{C}^*_{\mbqpNvars + 1} \cap \mathcal{NULL}}.
\end{mini}
As a quick aside, rewriting the problem in this way does not change the Lagrangian dual problem.
To see this, we first write the Lagrangian dual of Problem \eqref{eq:CPP_2cone} as
\begin{maxi}
{\gamma \in \mathbb{R}, \tilde{M} \in \mathbb{R}^{(\mbqpNvars + 1) \times (\mbqpNvars)}}{\gamma}
{}{}
\addConstraint{\tilde{M} = \quadmat{Q}{c}{\cdot} - \gamma \Quadmat{\cdot}{\cdot}{1} }
\addConstraint{\tilde{M} \in \mathcal{C}_{\mbqpNvars + 1} + \mathcal{NULL}^*}
\end{maxi}
and notice that $\mathcal{NULL}^*$ is spanned by 
\begin{equation}
    \left\{\quadmatrow{A^\top_{i, *} A_{i, *}}{-b_i A_{i, *}}{b_i^2}\right\} \cup \left\{\quadmat{-\basis{j} \basis{j}^\top}{\frac{1}{2} \basis{j}}{\cdot}]\right\}.
\end{equation}
Here, we take care to note that the \emph{Lagrangian} dual may differ from the \emph{conic} dual.
In particular, feasibility of \eqref{eq:MBQP} is sufficient for ensuring strong duality when $\tilde{M}$ is optimized over $(C^* \cap \mathcal{NULL})^*$ \cite[Prop 5.3.9]{bertsekas2009convex}; however, it is not guaranteed that $(C^* \cap \mathcal{NULL})^*$ is equal to $C + \mathcal{NULL}^*$; this is the case if and only if $C + \mathcal{NULL}^*$ is closed.

To establish strong duality, we will first assert that if Problem \eqref{eq:MBQP} is feasible with a bounded feasible region, then Problem \eqref{eq:CPP_hom} has a nonempty and bounded set of optimal solutions.
This follows directly from \cite[Corollary 2.6]{Burer2009}, which states that for all optimal solutions, $(x^*, X^*)$, of \eqref{eq:CPP_hom}, $x^*$ must lie within the convex hull of optimal solutions for \eqref{eq:MBQP}.
If the set of optimal solutions for \eqref{eq:MBQP} is non-empty and bounded, so is their convex hull, proving the boundedness of $x^*$.
Moreover, because the optimal solution may be factored in the form \begin{equation}\quadmat{X^*}{x^*}{1} = 
\sum_k
\begin{pmatrix} x^{(k)} \\ \xi^{(k)} \end{pmatrix}
\begin{pmatrix} x^{(k), \top} & \xi^{(k)} \end{pmatrix}
\end{equation}
with $\xi^{(k)} > 0$ and $\sum_k (\xi^{(k)})^2 = 1$ (by definition of copositivity), $X^*$ can be expressed as the sum of the outer products of optimal solutions of \eqref{eq:MBQP} with themselves.
In other words, because each $x^{(k)}$ is bounded, $X^* = \sum_k x^{(k)}x^{(k), \top}$ must be bounded as well.
This establishes that the set of optimal solutions of \eqref{eq:CPP_hom} is also non-empty and bounded, allowing us to apply \cite[Theorem 1.1]{kim2021strong}, which establishes strong duality of conic optimization problems with two cone constraints and a single hyperplane constraint under non-emptiness and boundedness of the optimal solution set.
The establishes that
\begin{equation}
    \max \eqref{eq:COP_hom} = \min \eqref{eq:CPP_hom}.
\end{equation}
While the optimal objective of \eqref{eq:COP_hom} may not be exactly attainable, there exists a feasible solution with objective value $\max \eqref{eq:COP_hom} - \epsilon$ where $\epsilon > 0$ can be arbitrarily small \cite{cifuentes2023sensitivity}.
This is not restrictive, as most numerical solvers only compute optimums with finite precision.
\end{proof}
\end{theorem}

While Theorem \ref{thm:hom_strong_duality} establishes strong duality of the homogenous CPP, we have yet to show that the non-homogeneous form also exhibits strong duality.
In order to do so, we will show that the supremum of the (non-homogeneous) copositive program upper-bounds that of the homogeneous program.
\begin{theorem}[Inhomogeneous Lower Bound]
The optimal objective of Problem \eqref{eq:COP} is at least that of Problem~\eqref{eq:COP_hom} (i.e., $\max \eqref{eq:COP} \geq \max \eqref{eq:COP_hom}$).
\begin{proof}
 We will do this by showing that for each $(\hat{\mu}, \hat{\lambda}, \hat{\gamma})$ there exists $(\mu, \lambda, \gamma)$ such that 
\begin{equation}
    M(\mu, \lambda, \gamma) = \hat{M}(\hat{\mu}, \hat{\lambda}, \hat{\gamma}),
\end{equation}
and 
\begin{equation}
    \gamma + \sum_i \mu^{\text{(lin)}}_i b_i + \mu^{\text{(quad)}}_i b^2_i = \hat{\gamma}.
\end{equation}
In other words, any feasible solution for \eqref{eq:COP_hom} can be transformed into a feasible solution for \eqref{eq:COP} with equal objective value.
To see this, we will suggestively break up $\gamma = \gamma^{\text{(res)}} + \sum_i \gamma_i$ so equation \eqref{def:M} can be expanded as
\begin{equation}
    \begin{aligned}
        M(\mu, \lambda, \gamma) = &\quadmat{Q}{c}{\cdot}\\
        & - \sum_i \left(\mu^{\text{(lin)}}_i \quadmatrow{\cdot}{\frac{1}{2} A_{i, *}}{\cdot} + \mu^{\text{(quad)}}_i \Quadmat{A^\top_{i, *} A_{i, *}}{\cdot}{\cdot} + \gamma_i \Quadmat{\cdot}{\cdot}{1} \right) \\
        & - \sum_{j \in B} \lambda_j \quadmat{-\basis{j} \basis{j}^\top}{\frac{1}{2} \basis{j}}{\cdot} - \gamma^{\text{(res)}} \Quadmat{\cdot}{\cdot}{1}
    \end{aligned}
\end{equation}
Then, the proposed $(\mu, \lambda, \gamma)$ is given by
\begin{align}
    \lambda_j &= \hat{\lambda}_j\\
    \mu_i^{\text{(lin)}} &= -2 b_i \hat{\mu}_i\\
    \mu_i^{\text{(quad)}} &= \hat{\mu}_i\\
    \gamma_i &= b_i^2 \hat{\mu}_i\\
    \gamma^{\text{(res)}} &= \hat{\gamma} 
\end{align}
Then, notice that 
\begin{align}
    &\mu^{\text{(lin)}}_i \quadmatrow{\cdot}{\frac{1}{2} A_{i, *}}{\cdot} + \mu^{\text{(quad)}}_i \Quadmat{A^\top_{i, *} A_{i, *}}{\cdot}{\cdot} + \gamma_i \Quadmat{\cdot}{\cdot}{1}\\
    &= -2 b_i \hat{\mu}_i\quadmatrow{\cdot}{\frac{1}{2} A_{i, *}}{\cdot} + \hat{\mu}_i \Quadmat{A^\top_{i, *} A_{i, *}}{\cdot}{\cdot} + b_i^2 \hat{\mu}_i \Quadmat{\cdot}{\cdot}{1}\\
    &= \hat{\mu}_i \left( -2 b_i \quadmatrow{\cdot}{\frac{1}{2} A_{i, *}}{\cdot} + \Quadmat{A^\top_{i, *} A_{i, *}}{\cdot}{\cdot} + b_i^2 \Quadmat{\cdot}{\cdot}{1}\right)\\
    &= \hat{\mu}_i \quadmatrow{A^\top_{i, *} A_{i, *}}{-b_i A_{i, *}}{b_i^2}
\end{align}
so by matching up terms in the sums, we see that $M(\mu, \lambda, \gamma) = \hat{M}(\hat{\mu}, \hat{\lambda}, \hat{\gamma})$.
As for the objective value, 
\begin{align}
    &\gamma + \sum_i \mu^{\text{(lin)}}_i b_i + \mu^{\text{(quad)}}_i b_i^2\\
    &= \gamma^{\text{(res)}} + \sum_i \mu^{\text{(lin)}}_i b_i + \mu^{\text{(quad)}}_i b_i^2 + \gamma_i\\
    &= \hat{\gamma}  + \sum_i -2 b^2_i \hat{\mu}_i +  b_i^2 \hat{\mu}_i + b_i^2 \hat{\mu}_i\\
    &= \hat{\gamma}  + \sum_i \hat{\mu}_i (-2 b^2_i  +  b_i^2  + b_i^2 )\\
    &= \hat{\gamma}
\end{align}
so for each $(\hat{\mu}, \hat{\lambda}, \hat{\gamma})$ the proposed $(\mu, \lambda, \gamma)$ has equal objective value.
\end{proof}
\end{theorem}

\begin{corollary}
If Problem \eqref{eq:MBQP} is feasible with bounded feasible region, then strong duality holds between Problems \eqref{eq:CPP} and \eqref{eq:COP}.
\begin{proof}
Theorem shows that $\max \eqref{eq:COP} \geq \max \eqref{eq:COP_hom}$.
 So we have $\max \eqref{eq:COP_hom} \leq \max \eqref{eq:COP} \leq \min \eqref{eq:CPP} = \min \eqref{eq:CPP_hom}$.
Combining this with $\max \eqref{eq:COP_hom} = \min \eqref{eq:CPP_hom}$ we get 
\begin{equation}
    \max \eqref{eq:COP_hom} = \max \eqref{eq:COP} = \min \eqref{eq:CPP} = \min \eqref{eq:CPP_hom}.
\end{equation}
Thus, strong duality must hold between \eqref{eq:COP} and \eqref{eq:CPP}.
\end{proof}
\end{corollary}

\subsection{Discretizing the copositivity checks}\label{subsec:disc_cop}
\subsubsection{Constructing the QUBO}
In this section, we will discuss forming the QUBO to approximate the copositivity checks. 
Formally, instead of solving \eqref{eq:norm_cop} with feasible region $\{z \in \mathbb{R}^{\mbqpNvars + 1}_{\scriptscriptstyle \geq 0} \mid \norm{z}_\infty \leq 1\}$, we will approximate the feasible region with $\{0, \frac{1}{\ndisc}, \ldots, \frac{\ndisc- 1}{\ndisc}, 1\}^{\mbqpNvars + 1}$, leading to a quadratic unconstrained integer optimization, 
\begin{mini}
{z}{z^\top M(\mu, \lambda, \gamma) z}
{\label{eq:disc_cop}\tag{QUIO}}{}
\addConstraint{z \in \left\{0, \frac{1}{\ndisc}, \ldots, \frac{\ndisc- 1}{\ndisc}, 1\right\}^{\mbqpNvars + 1}}
\end{mini}
For simplicity, assume that $\ndisc = 2^{\ndiscbin} - 1$ for some $\ndiscbin \in \mathbb{Z}_{\scriptscriptstyle > 0}$.
Then Problem \eqref{eq:disc_cop} is equivalent to minimizing \eqref{eq:qubo}, where
\begin{align}
\hat{M}(\mu, \lambda, \gamma) &= \Dmat^\top M(\mu, \lambda, \gamma) \Dmat
\end{align}
and
\begin{align}\label{eq:disc_coeff}
    \Dmat &:= \dfrac{1}{\ndisc}
    \begin{pmatrix}
    2^0 & \cdots & 2^{\ndiscbin - 1} &
    0 & \cdots & 0 &
    \cdots & 
    0 & \cdots & 0\\
    0 & \cdots & 0 &
    2^0 & \cdots & 2^{\ndiscbin  - 1} &
    \cdots &
    0 & \cdots & 0\\
    \vdots & \vdots & \vdots &
    \vdots & \vdots & \vdots &
    \vdots &
    \vdots & \vdots & \vdots\\
    0 & \cdots & 0 & 
    0 & \cdots & 0 &
    \cdots &
    2^0 & \cdots & 2^{\ndiscbin  - 1} 
    \end{pmatrix},
\end{align}
over the variables $\hat{z} \in \{0, 1\}^{k(n + 1)}$.
The new variables $\hat{z}$ simply represent the binary expansion of $z$, i.e., $z = \frac{1}{K} \Dmat \hat{z}$.
One could also use a unary expansion at the expense of a larger size expansion and redundancy in the encoding.
Additionally, while we have written out a uniform expansion for all variables, it is possible to have a heterogenous discretization scheme.
More sophisticated discretization schemes are discussed in depth in \cite{karimi2019practical}.

\subsubsection{Choosing a discretization size}
When discretizing the copositivity checks, it is critical to ensure that the discretization size is fine enough.
This section provides guidance for choosing a discretization size given a particular QUBO.
To do so, we will consider minimizing $z^\top M z$ over a discrete grid $z \in \{0, \frac{1}{\ndisc}, \ldots, \frac{\ndisc- 1}{\ndisc}, 1\}^n$ and bound the difference if we have minimized $\hat{z}^\top M \hat{z}$ over the hypercube instead, $\hat{z} \in [0, 1]^n$.
In particular, we are interested in the case where there are no certificates of non-copositivity on the discrete grid, i.e., $z^\top M z \geq 0$ for all $z \in \{0, \frac{1}{\ndisc}, \ldots, \frac{\ndisc- 1}{\ndisc}, 1\}^n$, yet there is $\hat{z} \in [0, 1]^n$ with $\hat{z}^\top M \hat{z} = -\delta < 0$.
We will decompose $\hat{z} = z + \Delta$ as the nearest grid point, $z ~\in \{0, \frac{1}{\ndisc}, \ldots, \frac{\ndisc- 1}{\ndisc}, 1\}^n$, plus a small correction factor $\Delta \in \mathbb{R}^n$.

Because the norm of this correction factor is bounded by $\norm{\Delta}_\infty \leq \frac{1}{2\ndisc}$ (i.e., by rounding), we will lower-bound $(z + \Delta)^\top M (z + \Delta)$ as a function of $\norm{\Delta}$.
Expanding $(z + \Delta)^\top M (z + \Delta)$ out we get 
\begin{equation}
    (z + \Delta)^\top M (z + \Delta) = z^\top Mz + 2 \Delta^\top M z + \Delta^\top M \Delta .
\end{equation}
Applying the assumption that $z^\top M z \geq 0$ for all $z \in \{0, \frac{1}{\ndisc}, \ldots, \frac{\ndisc- 1}{\ndisc}, 1\}^n$, we get
\begin{equation}
    (z + \Delta)^\top M (z + \Delta) \geq - |2 \Delta^\top M z + \Delta^\top M \Delta|.
\end{equation}
In order to lower-bound, $\hat{z}^\top M \hat{z}$, we upper-bound $|2 \Delta^\top M z + \Delta^\top M \Delta|$, 
\begin{align}
    | 2 \Delta^\top M z + \Delta^\top M \Delta| & = \norm{2 \Delta^\top M z + \Delta^\top M \Delta}_\infty \\
    &\leq 2 \norm{\Delta}_\infty \norm{z}_\infty \norm{M}_\infty +  \norm{\Delta}^2_\infty \norm{M}_\infty\\
    &\leq (2\norm{\Delta}_\infty +  \norm{\Delta}^2_\infty) \norm{M}_\infty\\
    &\leq \left(\frac{1}{K} + \frac{1}{4K^2}\right) \norm{M}_\infty.
\end{align}
Recall that we assume the minimum copositivity check is achieved at $-\delta = (z + \Delta)^\top M (z + \Delta)$, so if $\ndisc > \frac{1}{2(\sqrt{\frac{\delta}{\norm{M}_\infty} + 1} - 1)}$, then there exists $z \in \{0, \frac{1}{\ndisc}, \ldots, \frac{\ndisc- 1}{\ndisc}, 1\}^n$ with $z^\top M z < 0$.
This represents the coarsest discretization where optimizing over the discrete grid rather than the unit hypercube is insufficient for detecting the certificate of non-copositivity.

\subsection{Hyper-parameter optimization}
\label{subsec:hyperopt}
\begin{figure}
    \centering
    \includegraphics[width=\textwidth]{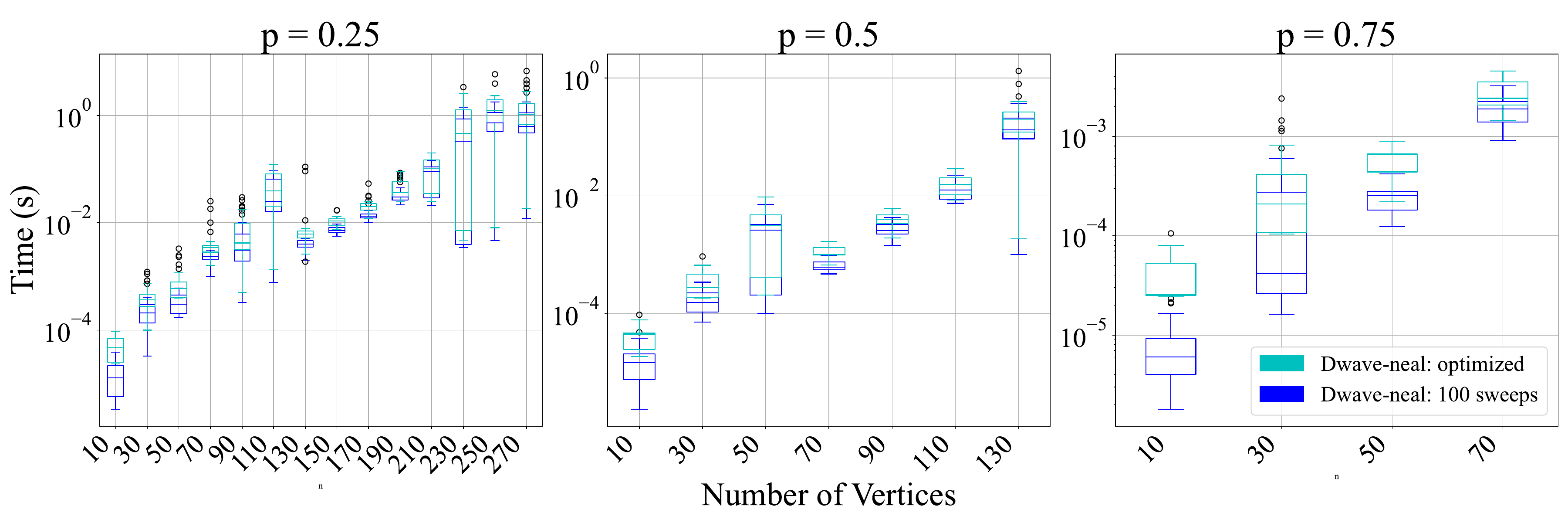}
    \caption{This figure plots the optimized $\texttt{TTT}_{0.999}$ when \texttt{Neal} was run with 100 sweeps. Optimization produces an order of magnitude speed-up for graphs with 10 nodes but does not result in significant speed-ups for larger graphs.}
    \label{fig:tuning_tts}
\end{figure}

\begin{figure}
    \centering
    \includegraphics[width=\textwidth]{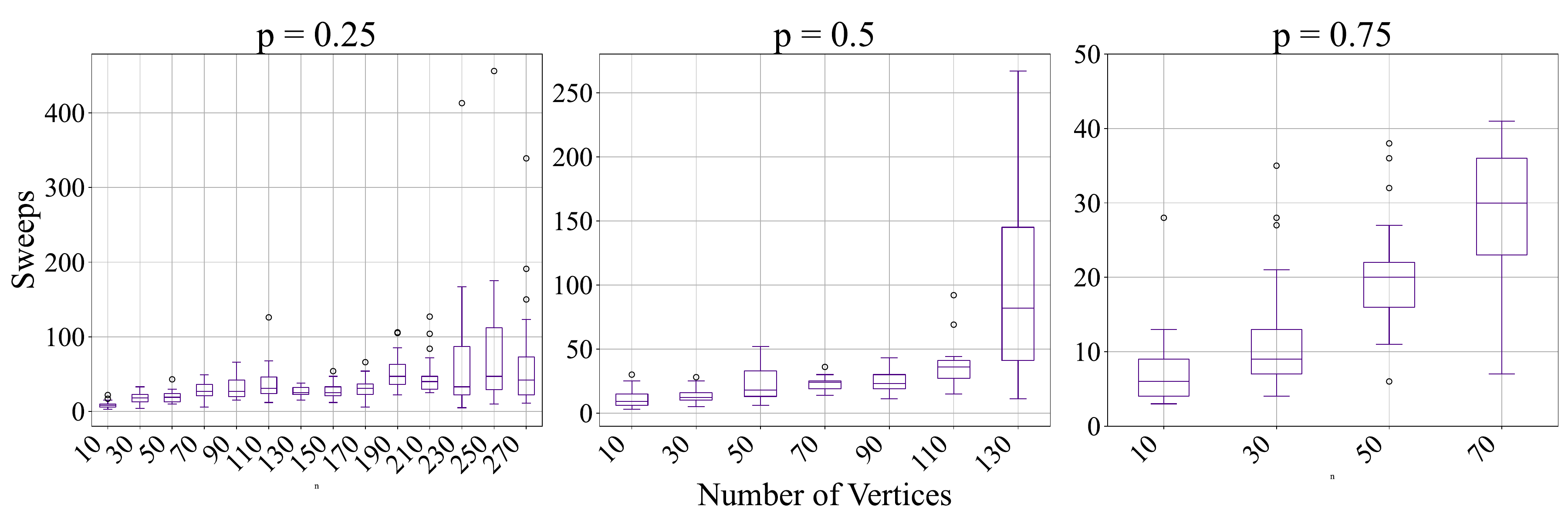}
    \caption{This figure plots the optimal number of sweeps for each of the problem instances.
    The optimal number of sweeps increases with the graph size. However, graphs of densities $p = 0.25$ and $p = 0.75$ require a comparable number of sweeps for graphs of the same size, while fewer sweeps are required for graphs with density $p = 0.5$}
    \label{fig:tuning_sweeps}
\end{figure}

To investigate further speed-ups from turning the simulated annealing parameters, we optimized the number of sweeps using \texttt{Hyperopt}~\cite{bergstra2013making} with 25 trials for each instance.
Figure \ref{fig:tuning_tts} plots the optimized $\texttt{TTT}_{0.99}$ and the $\texttt{TTT}_{0.99}$ when \texttt{Neal} was run with 100 sweeps.
While the optimization produced significant relative improvements for graphs with 10 nodes, the improvement for larger graphs remained marginal, especially in light of the computational overhead required to optimize the parameters.
Figure \ref{fig:tuning_sweeps} plots the optimal number of sweeps for each problem instance.
Generally, the optimal number of sweeps increases with the number of vertices.
While graphs of densities $p = 0.25$ and $p = 0.75$ require a comparable number of sweeps for graphs of the same size, fewer sweeps are required for graphs with density $p = 0.5$.

\subsection{Illustrative example}
In this section, we will walk through a small MBQP to illustrate the translation into the equivalent copositive program.
Consider the following mixed-binary optimization problem:
\begin{mini}
{x_1, x_2}
{
\begin{pmatrix} x_1 & x_2 \end{pmatrix}
\begin{pmatrix} 1 & -1\\ -1 & \cdot  \end{pmatrix}
\begin{pmatrix} x_1 \\ x_2 \end{pmatrix}
}
{\label{eq:ex_mbqp}\tag{Ex-MBQP}}{}
\addConstraint{x_1 + x_2 = 1}
\addConstraint{x_1, \, x_2 \in \mathbb{R}_{\scriptscriptstyle \geq 0}} .
\end{mini}
The optimal solution is given by 
\begin{equation}
    x_1^* = \frac{1}{3}, \quad x_2^* = \frac{2}{3},
\end{equation}
which gives an optimal objective value of $-\frac{1}{3}$.

\begin{figure}
\centering
\begin{tabular}{ccc}
  \includegraphics[clip, trim=8cm 0cm 7cm 1cm, width=0.25\textwidth]{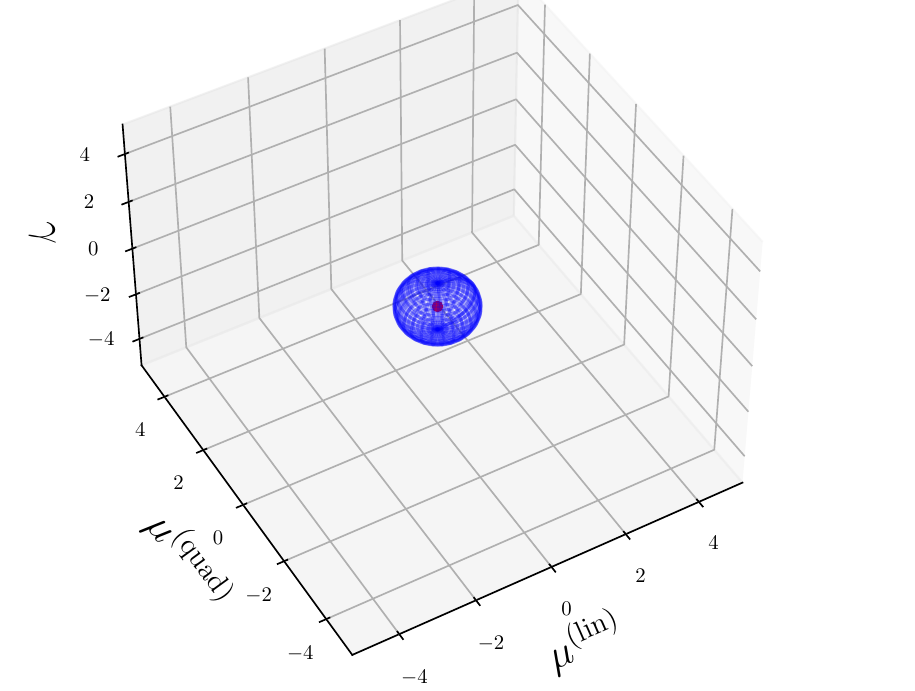}
  &
  \includegraphics[clip, trim=8cm 0cm 7cm 1cm, width=0.25\textwidth]{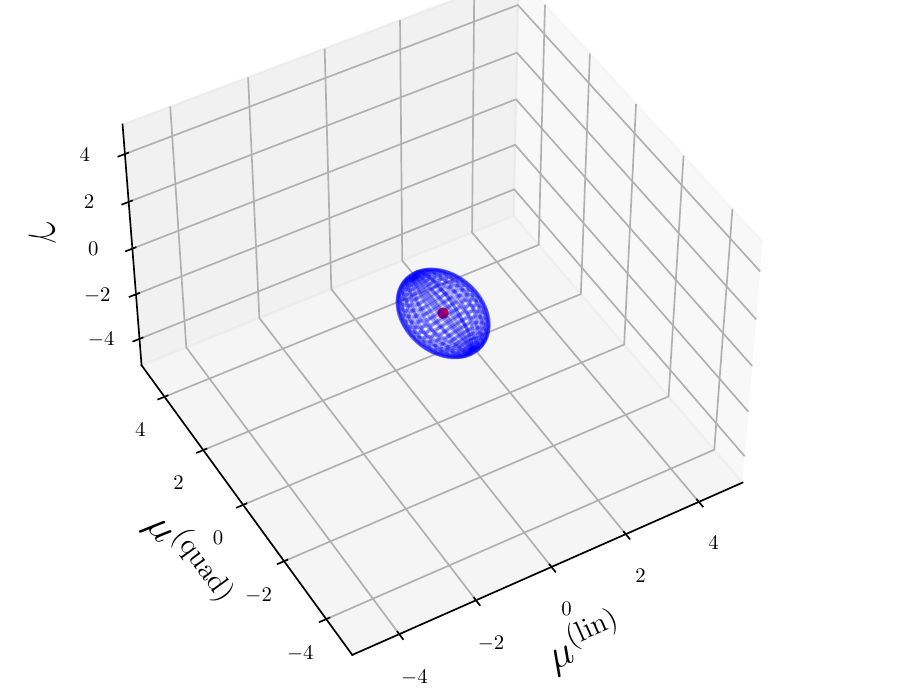} 
  &
  \includegraphics[clip, trim=8cm 0cm 7cm 1cm, width=0.25\textwidth]{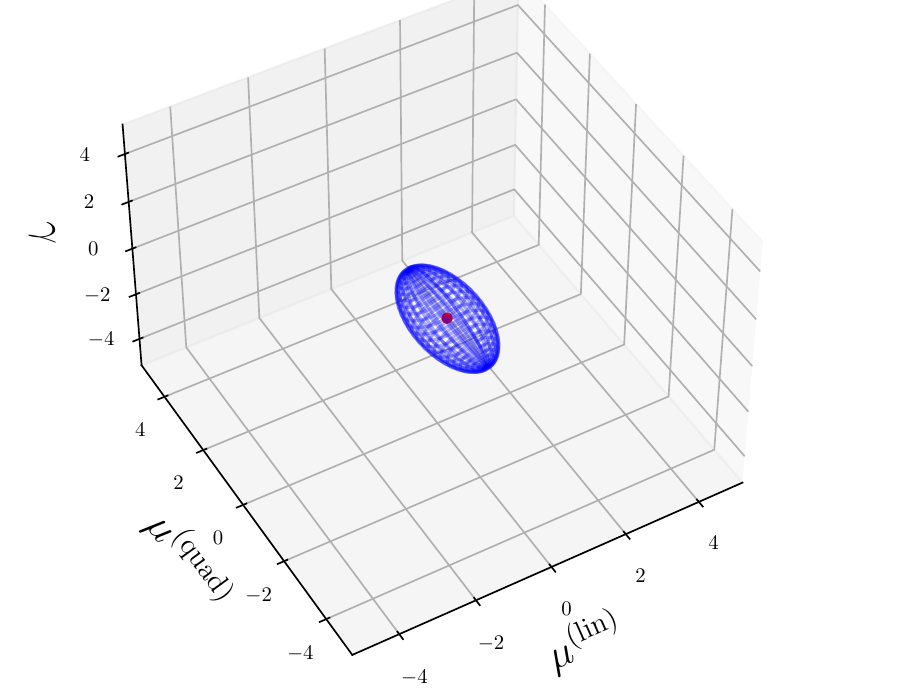}\\
 Iteration 0  & Iteration 1  & Iteration 2 \\
 \includegraphics[clip, trim=8cm 0cm 7cm 1cm, width=0.25\textwidth]{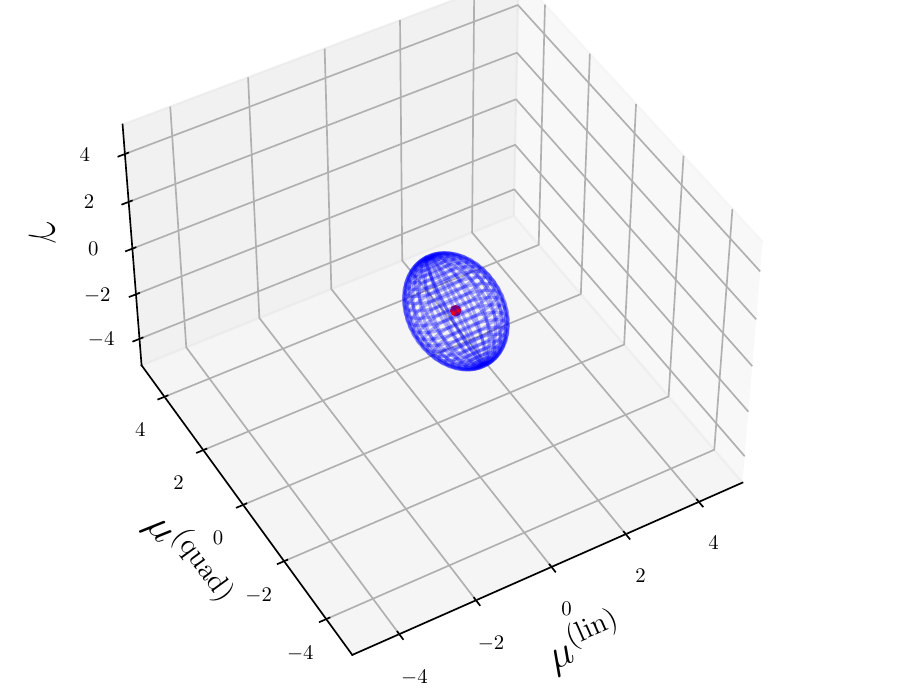}
  &
  \includegraphics[clip, trim=8cm 0cm 7cm 1cm, width=0.25\textwidth]{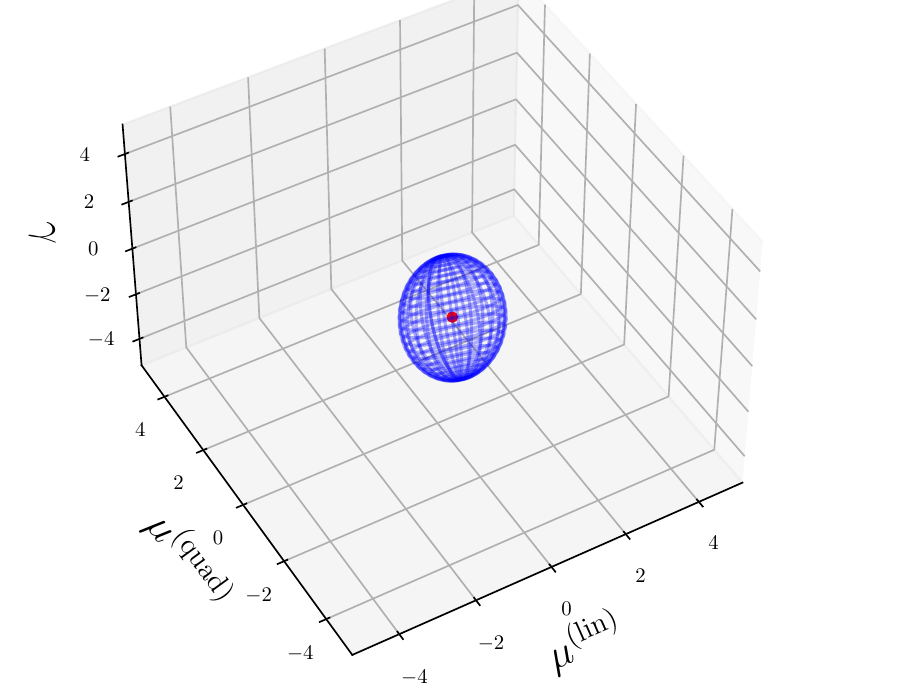} 
  &
  \includegraphics[clip, trim=8cm 0cm 7cm 1cm, width=0.25\textwidth]{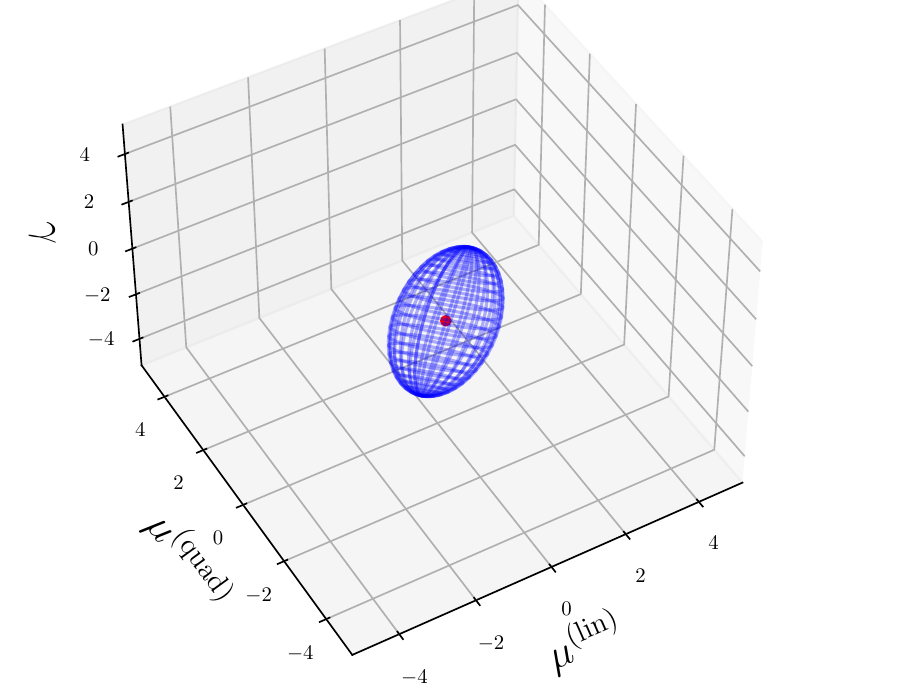}\\
 Iteration 3  & Iteration 4  & Iteration 5\\
 \includegraphics[clip, trim=8cm 0cm 7cm 1cm, width=0.25\textwidth]{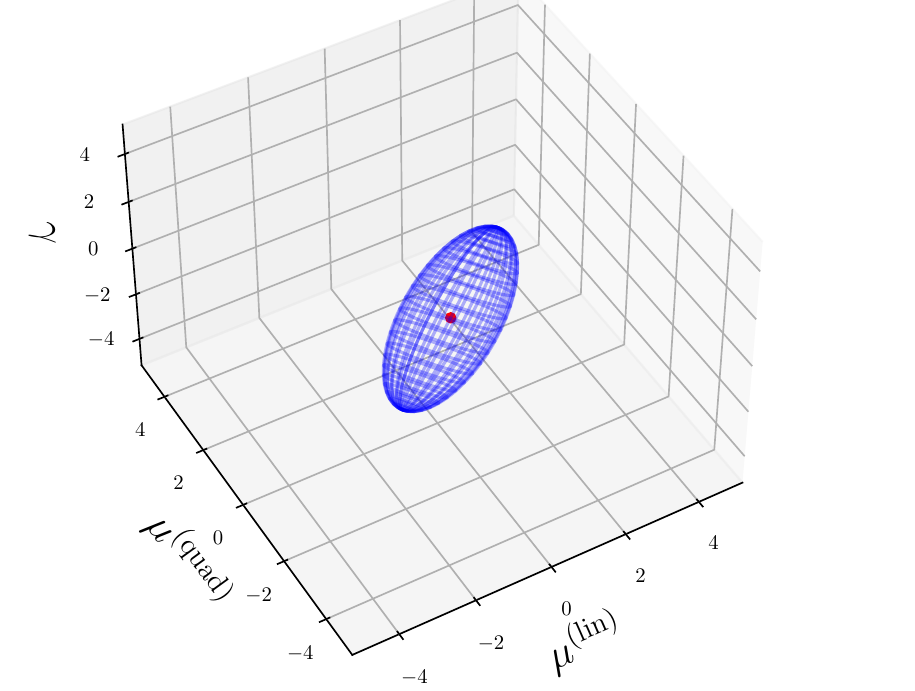}
  &
  \includegraphics[clip, trim=8cm 0cm 7cm 1cm, width=0.25\textwidth]{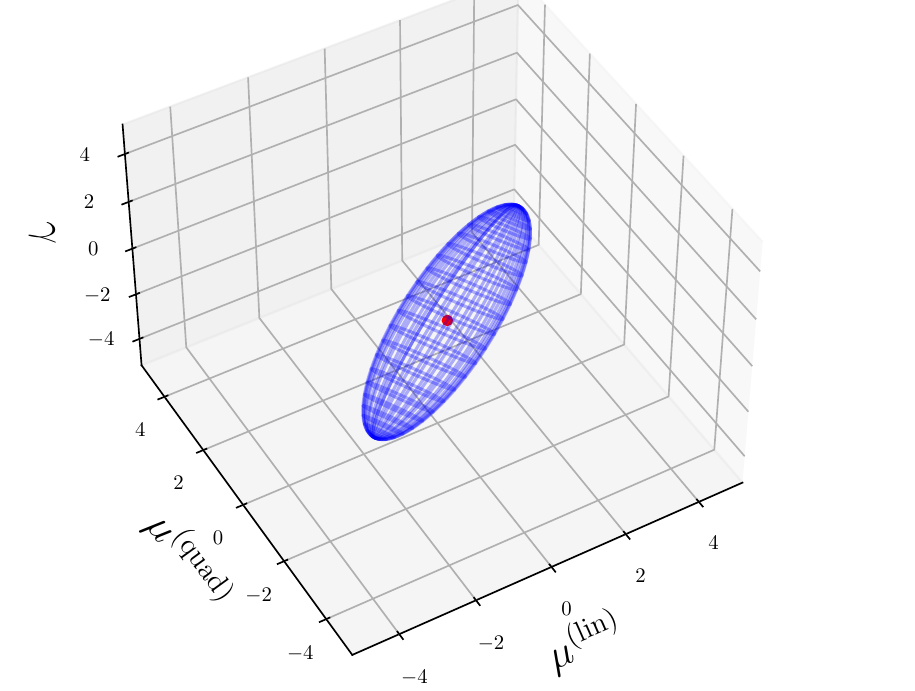} 
  &
  \includegraphics[clip, trim=8cm 0cm 7cm 1cm, width=0.25\textwidth]{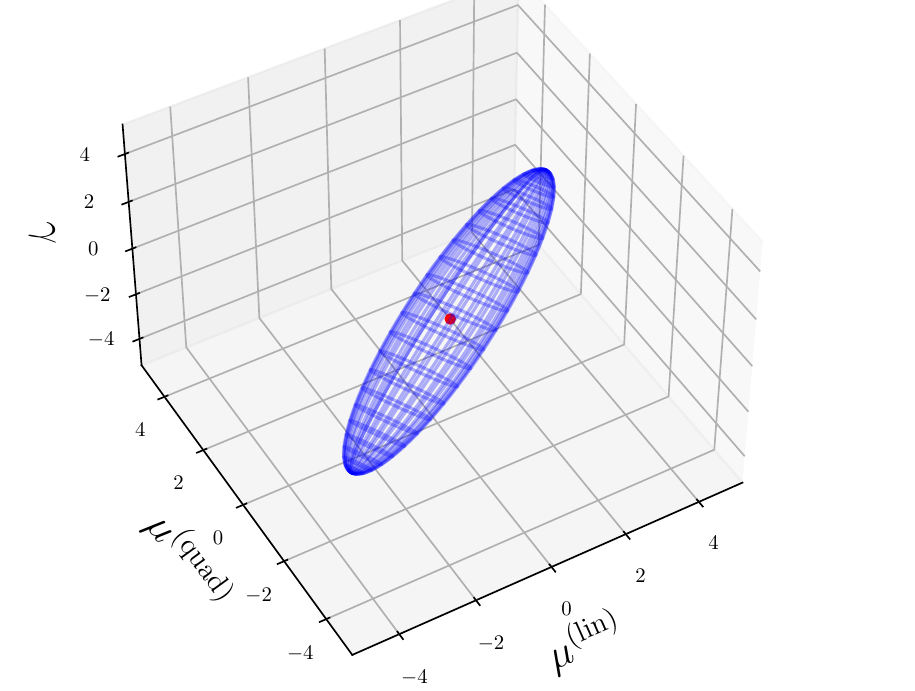}\\
 Iteration 6  & Iteration 7  & Iteration 8 
\end{tabular}
\caption{This figure plots the outer bounding ellipsoid for the first nine iterations of the copositive cutting plane algorithm applied to Problem \eqref{eq:ex_cop}.
In each plot, the red dot depicts the test point, and the blue ellipsoid plots the outer bounding ellipsoid at the start of the iteration.}
\label{fig:ex_ellipsoid}
\end{figure}

The equivalent completely positive program is given by 
\begin{mini}
{X \in \mathbb{R}^{2 \times 2}, x \in \mathbb{R}^2}
{
\innerMat{
\small\begin{pmatrix}
1 & -1 & \cdot \\ 
-1 & \cdot & \cdot \\ 
\cdot & \cdot & \cdot 
\end{pmatrix}
}{\quadmat{X}{x}{1}}
}
{\label{eq:ex_cpp}\tag{Ex-CPP}}{}
\addConstraint{
\innerMat
{\small\begin{pmatrix} 
\cdot&\cdot&1\\
\cdot&\cdot&1\\
1&1&\cdot\\
\end{pmatrix}}{\quadmat{X}{x}{1}} = 2}
\addConstraint{
\innerMat
{\small\begin{pmatrix} 
1&1&\cdot\\
1&1&\cdot\\
\cdot&\cdot&\cdot
\end{pmatrix}}{\quadmat{X}{x}{1}} = 1}
\addConstraint{
\innerMat
{\small\begin{pmatrix} 
\cdot&\cdot&\cdot\\
\cdot&\cdot&\cdot\\
\cdot&\cdot&1\\
\end{pmatrix}}{\quadmat{X}{x}{1}} = 1}
\addConstraint{\quadmat{X}{x}{1}\in \mathcal{C}^*_3.}
\end{mini}
The optimal solution of \eqref{eq:ex_cpp} is determined by the quadratic expansion of the optimal solution of \eqref{eq:ex_mbqp} as follows:
\begin{equation}
    \quadmat{X^*}{x^*}{1} = 
    \begin{pmatrix} x_1^* \\ x_2^* \\ 1 \end{pmatrix}
    \begin{pmatrix} x_1^* & x_2^*& 1 \end{pmatrix}
    =
    \begin{pmatrix}
    1/9 & 2/9 & 1/3\\
    2/9 & 4/9 & 2/3\\
    1/3 & 2/3 & 1
    \end{pmatrix}
\end{equation}
The dual copositive program is given by 
\begin{maxi}
{\mu, \lambda, \gamma}{\gamma + 2\mu^{\text{(lin)}} + \mu^{\text{(quad)}}}
{\label{eq:ex_cop}\tag{Ex-COP}}{}
\addConstraint{ M(\mu, \lambda, \gamma) \in \mathcal{C}_3}
\end{maxi}
where 
\begin{equation}
    M(\mu, \lambda, \gamma) = \begin{pmatrix}
1 -\mu^{\text{(quad)}} & -1-\mu^{\text{(quad)}} & -\mu^{\text{(lin)}} \\
-1-\mu^{\text{(quad)}} & -\mu^{\text{(quad)}} & -\mu^{\text{(lin)}} \\ 
-\mu^{\text{(lin)}} & -\mu^{\text{(lin)}} & -\gamma
\end{pmatrix}
\end{equation}
Figure \ref{fig:ex_ellipsoid} plots the outer bounding ellipsoid for the first 9 iterations of the copositive cutting plane algorithm (with the ellipsoid method as the cutting-plane algorithm) applied to Problem \eqref{eq:ex_cop}.
For each iteration, the red dot depicts the test point, and the blue ellipsoid plots the outer bounding ellipsoid at the start of the iteration.
The initial ellipsoid is chosen to be a sphere, but as the algorithm progresses, we observe that the outer bounding ellipsoid becomes elongated.
This behavior is explained by the fact that the optimal solution set for this particular problem is a line.

\subsection{Ellipsoid Algorithm}
In this section, we will overview the ellipsoid algorithm as a representative example of cutting-plane algorithms \cite{boyd2007localization}.
As suggested by its name, the outer approximation defined by an ellipsoid.
An ellipsoid in $\mathbb{R}^{\cpDims}$ is parametrized by a center, $x\in \mathbb{R}^{\cpDims}$, and positive definite matrix, $P \in \mathcal{S}^{\cpDims}_{++}$, and is defined as 
\begin{equation}
    \mathcal{E}(x, P) := \{s \in \mathbb{R}^{\cpDims} \mid (s - x)P^{-1}(s - x) \leq 1\}.
\end{equation}
The volume of $\mathcal{E}(x, P)$ scales with the determinant of $P$, 
\begin{equation}\label{eq:ell_vol}
    \texttt{Vol}(\mathcal{E}(x, P)) = \frac{\pi^{\cpDims/2}}{\Gamma\left(\dfrac{\cpDims}{2} +1 \right)} \sqrt{\text{det}(P)}.
\end{equation}
In the ellipsoid algorithm, the center will always be the test point. 
Given a separating hyperplane for $x$ (recall that this is defined by a vector $\shpVec \in \mathbb{R}^{\cpDims}$ such that $\shpVec^\top s \leq \shpVec^\top x $ for all $s \in S$), the ellipsoid updates the outer approximation with the minimum volume ellipsoid containing both $\mathcal{E}(x, P)$ and $\{ s \in \mathbb{R}^{\cpDims} \mid \shpVec^\top s \leq \shpVec^\top x \}$.
Conveniently, this ellipsoid, $\mathcal{E}(\hat{x}, \hat{P})$, has a closed form representation with
\begin{align}
    \hat{x}(x, P, a) &= x - \frac{P a}{(\cpDims + 1)\sqrt{a^{\top} P a}}, \\
    \hat{P}(x, P, a) &= \frac{\cpDims^2}{\cpDims^2 - 1}\left(P - \frac{2 P aa^\top P}{(\cpDims + 1)a^\top P a}\right).
\end{align}
Now we are in a position to present the ellipsoid algorithm in the terminology of Section \ref{subsec:cp_alg}.
To initialize the algorithm, the user chooses $x_0$ and $P_0$ appropriately and chooses a final tolerance $r$.
The outer approximations for each iteration are maintained via $x_k$ and $P_k$, i.e., $S_k = \mathcal{E}(x_k, P_k)$.
Evaluating the center of $S_k$ involves simply returning the stored value for $x_k$, i.e., $\texttt{Center}(S_k) = x_k$
The initial volume is determined from Equation \eqref{eq:ell_vol} as $R = \texttt{Vol}(\mathcal{E}(x_0, P_0))$.
At each iteration, the separation oracle, $\texttt{Oracle}(x_k)$, is evaluated using Algorithm \ref{alg:sep_or}, and outer approximation is updated as follows:

\begin{align}
x_{k + 1} &= \hat{x}(x_k, P_k, \texttt{Oracle}(x_k)), \\
P_{k + 1} &= \hat{P}(x_k, P_k, \texttt{Oracle}(x_k)),\\
\texttt{Add\_Cut}(S_k, \texttt{Oracle}(x))
&= \mathcal{E}(x_{k + 1}, P_{k + 1}).
\end{align}

\subsection{Additional plots}

\begin{figure}[h!]
    \centering
    \begin{minipage}[b]{\linewidth}
    {
    \label{fig:min_max_gurobi_5}
    \includegraphics[width=\textwidth]{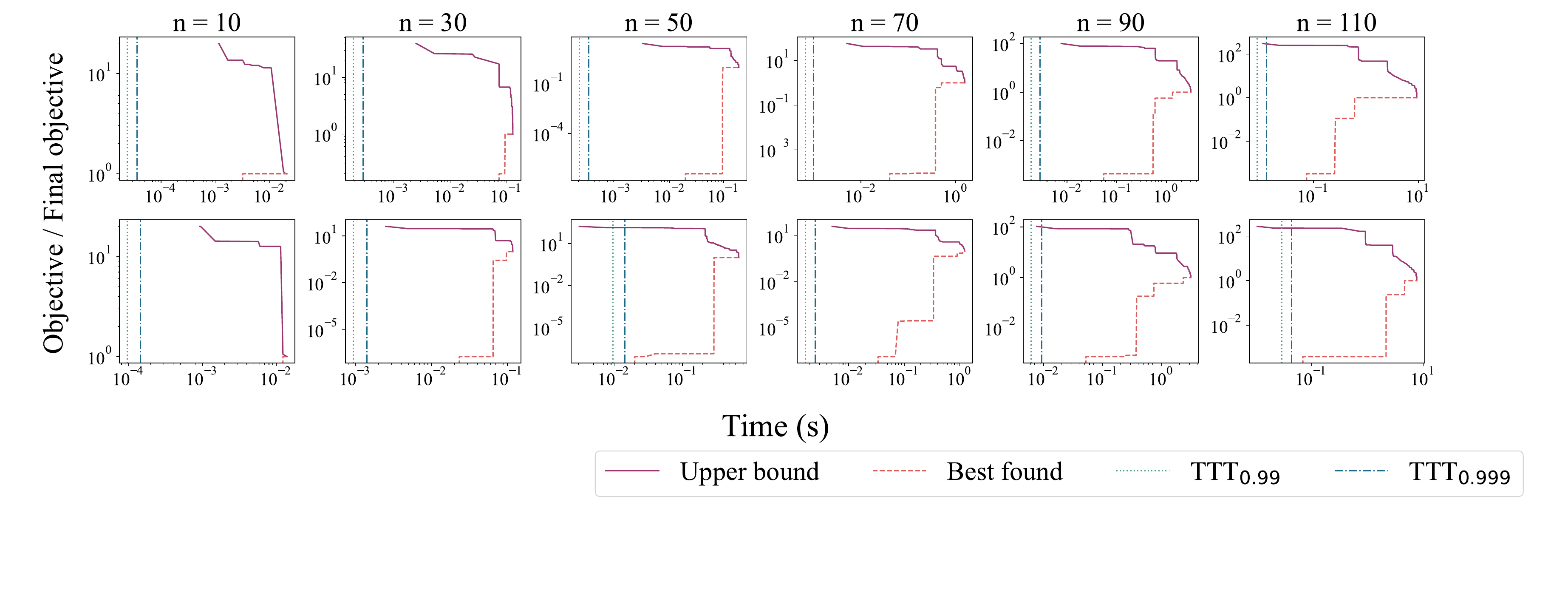}
    }
    \end{minipage}\par\vspace{-2\baselineskip}
    \centering
    \begin{minipage}[b]{\linewidth}
    {
    \label{fig:min_max_gurobi_75}
    \includegraphics[width=\textwidth]{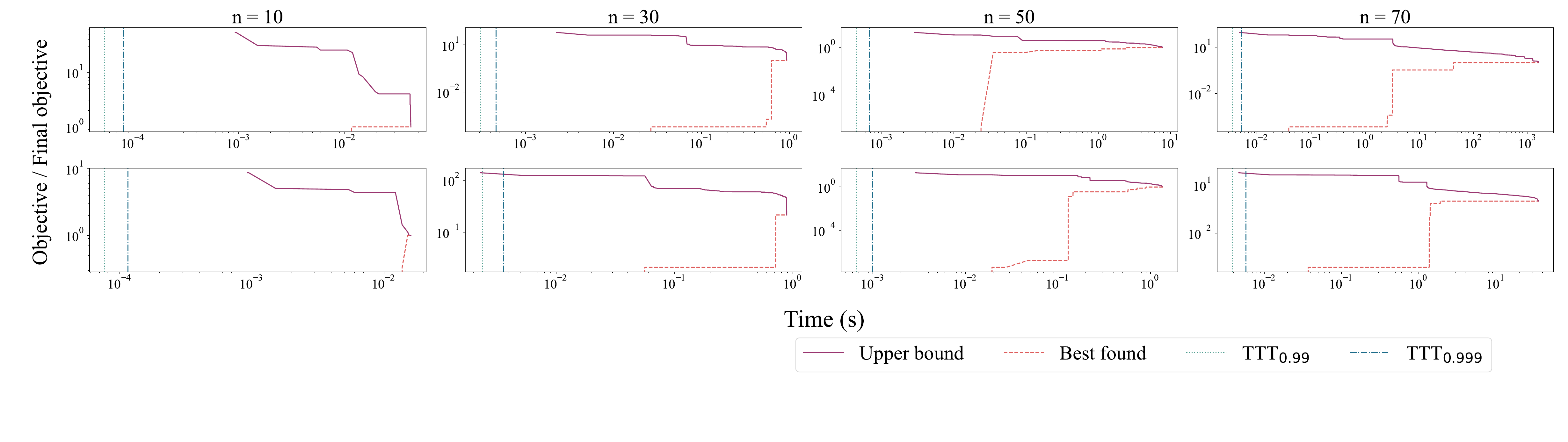}
    }
    \end{minipage}\par\vspace{-2\baselineskip}
    \caption{This figure depicts sample trajectories of \texttt{Gurobi}'s upper and lower bounds against $\texttt{TTT}_{0.99}$ and $\texttt{TTT}_{0.999}$ for edge density $p=0.5$ (above) and  $p=0.75$ (below).
    For each graph size, the top row represents the instance where the ratio between \texttt{Gurobi}'s solution time and $\texttt{TTT}_{0.99}$ is the greatest, and the bottom row represents the instance where the ratio is the smallest--all instances were run with 100 sweeps.
    In most instances, \texttt{Neal} reaches the $\texttt{TTT}_{0.999}$ confidence before \texttt{Gurobi} even returns a callback.}
\end{figure}

\begin{figure}[h!]
    \centering
    \includegraphics[width=\textwidth]{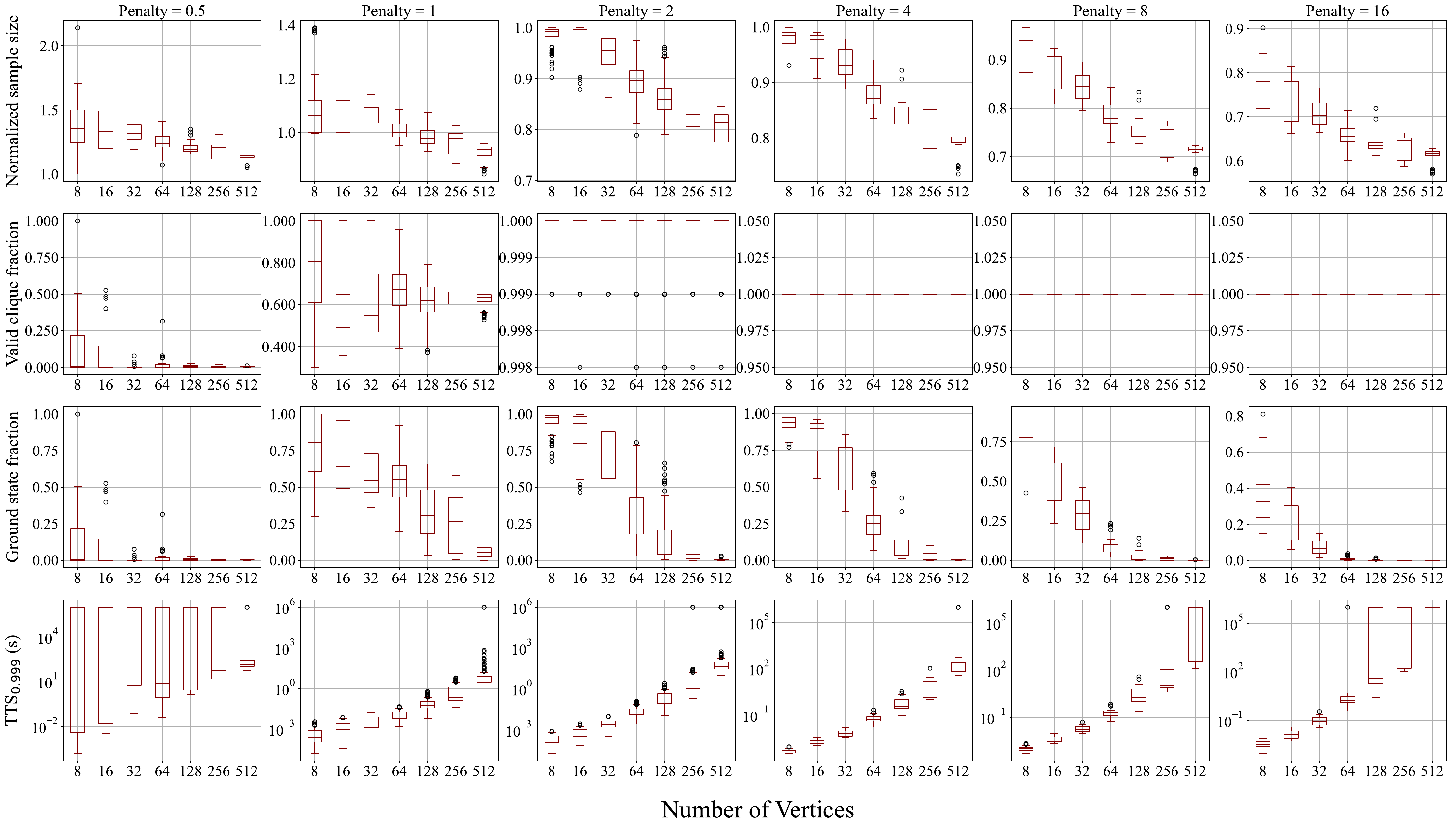}\par\vspace{-1\baselineskip}
        \caption{This figure plots the normalized sample size (the size of the returned solution divided by the ground truth maximum clique size) and the fraction of reads that resulted in a valid clique for graph density $p = 0.5$. These figures were used to compute the fraction of reads resulting in a ground state solution and the corresponding $\texttt{TTT}_{0.999}$ (also plotted). As the penalty weight is increased, the normalized sample size decreases, and the fraction of valid cliques increases. This highlights the delicate trade-off between constraints and the objective in penalty formulations.}
    \label{fig:penalty_05}
\end{figure}

\begin{figure}
    \centering
    \includegraphics[width=\textwidth]{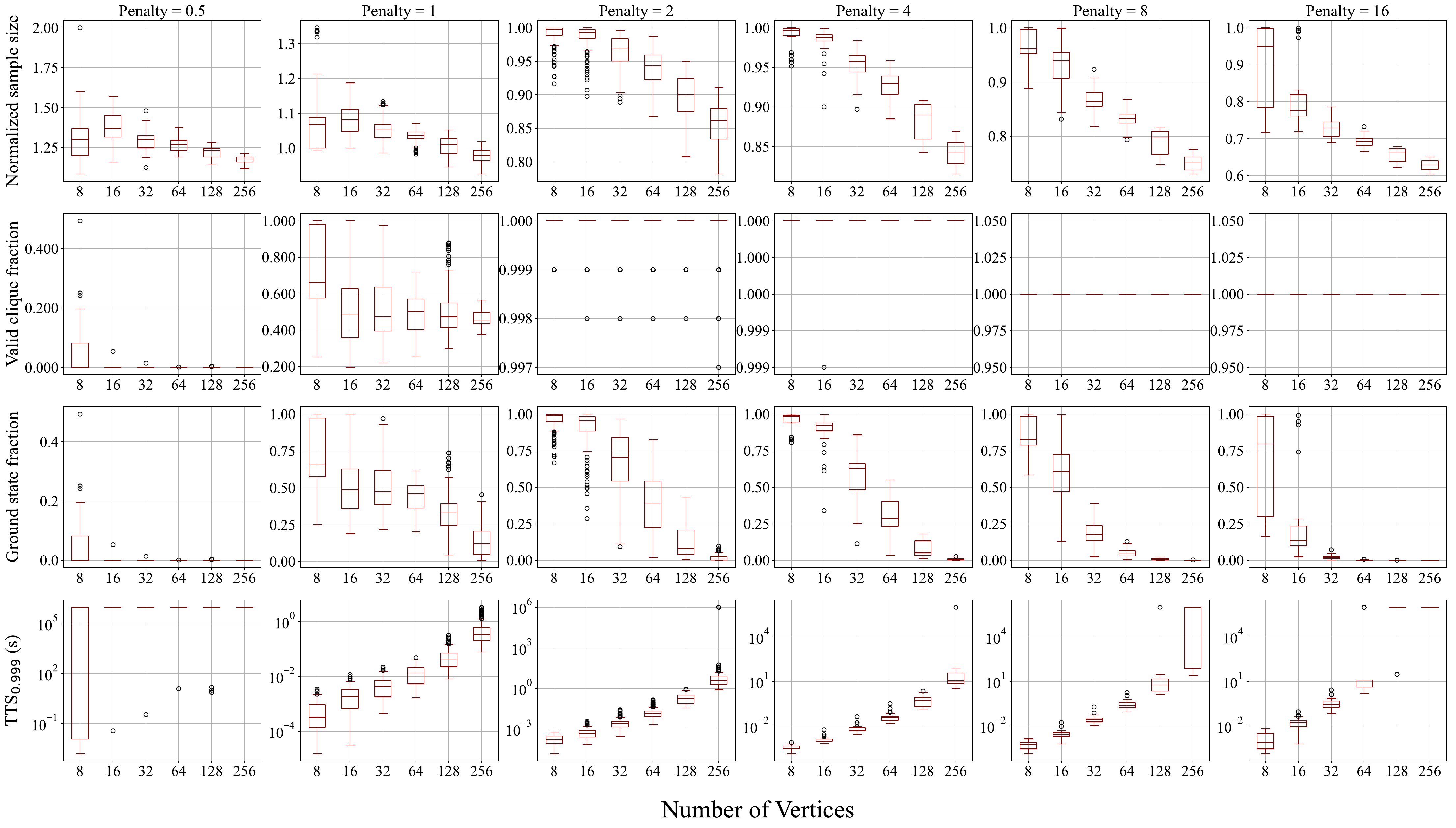}\par\vspace{-1\baselineskip}
        \caption{This figure plots the normalized sample size (the size of the returned solution divided by the ground truth maximum clique size) and the fraction of reads that resulted in a valid clique for graph density $p = 0.75$. These figures were used to compute the fraction of reads resulting in a ground state solution and the corresponding $\texttt{TTT}_{0.999}$ (also plotted). As the penalty weight is increased, the normalized sample size decreases, and the fraction of valid cliques increases. This highlights the delicate trade-off between constraints and the objective in penalty formulations.}
    \label{fig:penalty_75}
\end{figure}

\end{document}